\documentclass[12pt,leqno]{article}
\usepackage[a4paper, margin=24.0mm]{geometry}
\usepackage{amssymb}
\usepackage{amsmath}
\usepackage{amsthm}
\usepackage{stmaryrd}
\usepackage{mathtools,mathrsfs}
\usepackage[pagebackref]{hyperref}
\hypersetup{citecolor=purple, linkcolor=blue, colorlinks=true}
\usepackage{booktabs}
\usepackage{enumerate}
\usepackage{csquotes}
\usepackage{leftidx}
\usepackage{tikz}
\usepackage{url}
\usepackage{comment}
\usepackage[T2A,T1]{fontenc}
\usepackage{xfrac}
\usepackage{authblk}  

\newtheorem{theorem}{Theorem}[section]

\newtheorem*{proposition*}{Proposition}

\newtheorem{lemma}[theorem]{Lemma}
\newtheorem{corollary}[theorem]{Corollary}
\theoremstyle{definition}
\newtheorem{definition}[theorem]{Definition}
\newtheorem{remark}[theorem]{Remark}
\newtheorem{example}[theorem]{Example}

\newtheorem{hypothesis}[theorem]{Hypothesis}

\newcommand{\A}{{\mathcal A}}

\newcommand{\Aut}{\textup{Aut}}
\newcommand{\CAut}{\textup{CAut}}
\newcommand{\Char}{\textup{char}}
\newcommand{\Cyc}{\textup{C}}
\newcommand{\eps}{\varepsilon}
\newcommand{\Hol}{\textup{Hol}}
\newcommand{\Hom}{\textup{Hom}}
\newcommand{\GammaL}{\Gamma\textup{L}}

\newcommand{\GammaSp}{\Gamma\textup{Sp}}
\newcommand{\CSp}{\textup{CSp}}
\newcommand{\Gal}{\textup{Gal}}
\newcommand{\GL}{\textup{GL}}
\newcommand{\SL}{\textup{SL}}

\newcommand{\PSL}{\textup{PSL}}

\newcommand{\Sp}{\textup{Sp}}
\newcommand{\F}{\mathbb{F}}
\newcommand{\Fit}{\textup{F}}
\newcommand{\restrict}{\kern-2pt\downarrow\kern-2pt}

\newcommand{\Sym}{\textup{Sym}}
\newcommand{\Tr}{\textup{Tr}} 
\newcommand{\V}{{\mathcal V}}    
\newcommand{\W}{{\mathcal W}}    
\newcommand{\Z}{\textup{Z}}
\newcommand{\ZZ}{\mathbb{Z}}

\renewcommand{\ge}{\geqslant}
\renewcommand{\le}{\leqslant}

\newcommand{\lhdeq}{\trianglelefteqslant}    
\newcommand{\ord}[1]{\textup{ord}(#1)}
\newcommand{\order}{\textup{ord}}
\newcommand{\coloneq}{\vcentcolon=}

\begin{document}

\title{Classifying finite groups \texorpdfstring{$G$}{} with three \texorpdfstring{$\Aut(G)$}{}-orbits}

\author[1]{Stephen P. Glasby\footnote{Supported by the Australian Research Council (ARC) Discovery Project DP190100450}}
\affil[1]{\small Center for the Mathematics of Symmetry and Computation,
  University of Western Australia, Perth 6009, Australia.
  \href{mailto:Stephen.Glasby@uwa.edu.au}{Stephen.Glasby@uwa.edu.au}
	}


\maketitle

\abstract{
  We give a complete and irredundant list of the finite groups $G$
  for which $\Aut(G)$, acting naturally on $G$, has precisely $3$
  orbits. There are 7 infinite families: one abelian, one non-nilpotent, three
  families of non-abelian $2$-groups and two families of non-abelian
  $p$-groups with $p$ odd. The non-abelian $2$-group examples were first
  classified by Bors and Glasby in 2020 and non-abelian $p$-group examples
  with $p$ odd were classified independently by Li and Zhu~\cite{LZ2}, and by
  the author, in March 2024.
  \vskip2mm\noindent 
  {\bf Dedication: } To Otto H. Kegel on the occasion of his 90th birthday
  \vskip2mm\noindent 
}



\section{Introduction}\label{s:intro}

This work was presented at the Ischia Group Theory Conference
on 8 April~2024.

A group $G$ is called a $k$-orbit group if $\Aut(G)$, acting naturally on $G$,
has precisely $k$ orbits.
The cyclic group $\Cyc_{p^{k-1}}$, $k\ge1$, and the generalized quaternion
group $Q_{2^{k-1}}$ of order $2^{k-1}$, $k\ge5$, are both $k$-orbit $p$-groups.
Although we assume that $G$ is a \emph{finite} group, many of our examples
generalize to infinite groups. (For an infinite group $G$ it is sometimes
natural to count orbits under the subgroup of \emph{topological automorphisms}.)
Write $\ord{G}\coloneq\{|g|\mid g\in G\}$ for the set of
element orders in $G$. A 2-orbit group $G$ is elementary abelian and
$\ord{G}$ equals $\{1,p\}$ where $p$ is a prime.
If $G$ is a 3-orbit group, then $\ord{G}$ equals
$\{1,p\}$, with two orbits of elements of order $p$, or $\{1,p,p^2\}$
or $\{1,p,r\}$ where $p$ and $r$ are distinct primes. Thus 
a 3-orbit group is solvable by Burnside's $p^\alpha r^\beta$ theorem.
A complete $k$-orbit group $G$ has $\Aut(G)=\textup{Inn}(G)\cong G$ and
$k$ equals the number $k(G)$ of
$G$-conjugacy classes. In general, for a (finite) $k$-orbit group $G$ we have
the bounds $|\ord{G}|\le k\le \min\{k(G),|G|/p\}$ where $p$ is the
smallest prime divisor of $|G|$. The upper bound $|G|/p$ is due to~\cite[Theorem~1]{LM}. A second (obvious) lower bound is
  $1+\left\lceil(k(G)-1)/|\textup{Out}(G)|\right\rceil\le k$, see~\cite[Remark~4.1]{JKOB}.

This paper gives a complete and irredundant list of finite 3-orbit groups.
Clearly $\{1\}$ is an orbit, and for each characteristic subgroup $N$ of $G$ the
set $N^{\#} :=N\setminus\{1\}$ is a union of orbits.  If a $k$-orbit group
has a chain $G=N_1>\cdots>N_k=\{1\}$ of $k$ characteristic subgroups,
then the sets $N_i\setminus N_{i+1}$ for $1\le i<k$ and $\{1\}$ are the
$k$ non-trivial $\Aut(G)$-orbits, see Lemma~\ref{L:C} for more information.

Let $G$ be a 3-orbit group and set $N\coloneq\langle G',\Phi(G)\rangle$.
Then $N\ne1$ as $G$ is not elementary abelian, and $N\ne G$ as $G$ is solvable.
Hence $G/N$ and $N$ are both 2-orbit groups and thus elementary abelian.
We write $V=G/N=\F_r^{\,m}$ and $W=N=\F_p^{\,n}$ where $r,p$ are primes (possibly
equal). Observe that $\Aut(G)$ induces on $V$ a linear subgroup $A\le\GL_m(r)$,
and induces (via restriction) a linear subgroup $B\le\GL_n(p)$. Furthermore,
$A$ is transitive on the set $V\setminus\{0\}$ of non-zero vectors of $V$,
and $B$ is transitive on the set $W\setminus\{0\}$ of non-zero vectors of $W$.
We use Hering's theorem~\cite[\S5]{cH},~\cite[p.\,512]{Lie87a} to classify
the linear subgroups $A$ and $B$ (see Theorem~\ref{T:HeringNew}).

  We give a complete (Theorem~\ref{T:main}) and irredundant
  list (Theorem~\ref{T:irred}) of finite 3-orbit groups. Our list in
  Table~\ref{T:GVAWB} agrees with Theorem~B of Li and Zhu~\cite[Theorem~B]{LZ2}
as line 7 of Table~\ref{T:GVAWB} corresponds to their families (7) and (8).
Hering's theorem is a key tool in~\cite{BG,LZ1,LZ2} and the present
paper, our version is Theorem~\ref{T:HeringNew}. Consequently, the solvable residual plays an important role in
  our proof and in~\cite{LZ2}, as does the representation theory of the exterior square of a module.  Remark~\ref{R:irred} describes some similarities
  with our approach and that in~\cite{LZ2}. The substantial differences involve
  Theorem~\ref{T:irred} and a more geometric approach to the $p$-groups
  on line~7 of Table~\ref{T:GVAWB}
  below. We were also motivated by the possibility of classifying
4-orbit groups. We find it convenient to
represent our $p$-group (and solvable group) examples as Cartesian products
of vector spaces, with specified multiplication rules, and linear actions.

\begin{theorem}\label{T:main}
  Let $G$ be a finite $3$-orbit group with
  $N=\langle G',\Phi(G)\rangle$ and $|N|=p^n$.
  Then $1<N<G$ and $G$ is isomorphic to a group
  in lines $1-7$ of Table~$\ref{T:GVAWB}$. Moreover, the values
  of $V\cong G/N,A=\Aut(G)^V,W\cong N,B=\Aut(G)\restrict W$ are valid,
  where $\Aut(G)^V$ and $\Aut(G)\restrict W$ denote the groups
  induced on $G/N$ and $N$ by $\Aut(G)$.
\end{theorem}

\renewcommand{\arraystretch}{1.2}
\begin{table}[!ht]
\label{T:GVAWB}
\caption{3-orbit groups $G$ and $V\cong G/N,A=\Aut(G)^{G/N},B=\Aut(G)\restrict N$}
\begin{center}
\begin{tabular}{llclcll}\hline
  $\ \ \ \;G$&$V$&$A$&$N$&$B$&Comments&Ref.\\ \hline
  1. $(\Cyc_{p^2})^n$&$\F_p^{\;n}$&$\GL_n(p)$&$\F_p^{\;n}$&$\GL_n(p)$&$p\ge2$, $G$ abelian&p.\;\pageref{P:line1}\\
  2. $\F_{q}^{\;d}\rtimes\Cyc_r$&
  $\F_r$&$\GL_1(r)$&$\F_{q}^{\;d}$&$\GammaL_d(q)$&$q\kern-1.5pt=\kern-1.5ptp^{r-1}$\kern-4pt,\,$p\kern-1.5pt\ne\kern-1.5pt r$,\,$d\kern-1.5pt=\kern-1.5pt\frac{n}{r-1}$& \ref{E:Gpr}\\
  3. $A(n,\theta)$&$\F_2^{\;n}$&$\GammaL_1(2^n)$&$\F_2^{\;n}$&$\GammaL_1(2^n)$&Def.~\ref{D:ABP}(a), $n\kern-2.2pt\ne\kern-2.2pt 2^\ell$& \ref{L:abc}\\
  4. $B(n)$&$\F_2^{\;2n}$&$\GammaL_1(2^{2n})$&$\F_2^{\;n}$&$\GammaL_1(2^n)$&Def.~\ref{D:ABP}(b), $n\kern-1pt\ge\kern-1pt1$& \ref{L:abc}\\
  5. $P$&$\F_2^{\;6}$&$\Cyc_7\rtimes\Cyc_9$&$\F_2^{\;3}$&$\GammaL_1(2^3)$&Def.~\ref{D:ABP}(c), $n=3$& \ref{L:abc}\\
  6. $\F_{q}^{\;3}\kern-3pt:\kern-2pt\F_{q}^{\;3}$&$\F_{q}^{\;3}$&$\GammaL_3(q)$&$\F_{q}^{\;3}$&$\GammaL_3^{+}(q)$&$q=p^{\frac{n}{3}}$\ {\rm odd}, $3\mid n$& \ref{L:GL3}\\
  7. $\F_{p^n}\kern-3pt:\kern-2pt\F_{q}^{\;\frac{m}{b}}$&$\F_{q}^{\;\frac{m}{b}}$&$\Sp_{\frac{m}{b}}(\kern-1pt q\kern-1pt)\kern-2pt\le$&$\F_{p^n}$&$\GammaL_1(p^n\kern-1pt)\kern-2pt\le$&$q=p^b$\ {\rm odd}, $n\mid b\mid m$& \ref{L:SpTr}\\ [0.8mm]
  \hline
\end{tabular}
\end{center}
\end{table}

The number $\omega(G)$ of $\Aut(G)$-orbits on $G$ can sometimes be
calculated without knowing all of $\Aut(G)$, see Lemma~\ref{L:A}.
Nonsolvable $k$-orbit groups have been classified for $1\le k\le6$.
There are none when $k\le3$.
If $G$ is finite and $S$ is a composition factor of $G$, then $S$ is simple and
$\omega(G)\ge\omega(S)$. If $G$ is a finite nonsolvable $k$-orbit group,
then $G=A_5$ if $k=4$ by~\cite{LM},
$G\in\{\PSL_2(q)\mid q\in\{7,8,9\}\}$ if $k=5$ by~\cite{Stroppel}, and
$G\in\{\PSL_3(4),\F_4^2\rtimes\SL_2(4)\}$ if $k=6$ by~\cite{DGB}.
By contrast,
classifying \emph{solvable} $k$-orbit groups is extremely difficult. This
paper focuses on classifying (solvable) 3-orbit groups while contemplating
$4$-orbit groups, see Lemma~\ref{L:GL3} and Section~\ref{S:4orbit}.

We note that certain permutation groups of `rank'~$k$ give rise to
$k$-orbit groups.
A permutation group $G\le\Sym(\Omega)$ has \emph{rank} $k$ if it is
transitive on $\Omega$ and the point stabilizer $G_\omega$ 
has $k$ orbits including $\{\omega\}$; equivalently $G$ has $k$ orbits on
$\Omega\times\Omega$.

\begin{lemma}\label{L:rank}
  If $\Aut(G)$ has $k$ orbits on $G$, then the subgroup
  $\Hol(G)=G\rtimes\Aut(G)$ of the symmetric group $\textup{Sym}(G)$
  has rank $k$ and the stabilizer $\Hol(G)_1$ of $1$ is $\Aut(G)$.
\end{lemma}

Although rank~2 permutation groups have been classified, rank 3 groups
have only been classified in certain cases e.g.~\cite{Dornhoff,Foulser69}, or
when they are is quasiprimitive~\cite{DGLPP}, or
innately transitive~\cite{BDP}. Sadly for us,
the holomorph\footnote{A permutation group $H$ with a regular normal subgroup $G$ satisfies $G\lhdeq H\le\Hol(G)\le\Sym(G)$.} $\Hol(G)$ of a 3-orbit group $G$ 
is commonly not solvable, and $\Hol(G)$ is never innately transitive 
as $N=\langle \Phi(G),G'\rangle$ is intransitive and is the unique minimal normal subgroup
of $\Hol(G)\le\Sym(G)$. For a recent history of the classification of certain
low rank groups see~\cite[pp.\,177--178]{GMZ}. However, if $G$ is a 3-orbit $2$-group, then $\Aut(G)$
\emph{is} solvable by~\cite[Proposition~3.1]{BG}. Hence is $\Hol(G)$ also
solvable.
This is reminiscent of Theorem 2 of~\cite{Bryukhanova} by
Bryukhanova, and is used to prove Theorem~\ref{T:Ainftyeq1} when $p=2$.
The 3-orbit 2-groups $G$ have been classified
in~\cite[Corollary~1.3]{LZ1} and~\cite[Theorem~1.2]{BG}.
We find it useful to explicitly specify the (linear) action of $\Aut(G)$
on the vector space $N\cong\F_p^{\,n}$.

  We remark that Theorem~\ref{T:HeringNew} lists the perfect transitive
  linear groups and includes $\SL_2(5)\le\GL_4(3)$ which is not listed
  in~\cite[Theorem~3.1]{LZ2}. Fortunately, the example
  $\SL_2(5)\le 2_{-}^{1+4}.A_5\le\GL_4(3)$
  does not lead to an example of a 3-orbit group. The last sentence of
  Theorem~\ref{T:HeringNew} importantly follows from the converse
  in~\cite[Appendix~A]{Lie87a}.

Our proof of Theorem~\ref{T:main} is divided into three cases.
Section~\ref{S:prep} outlines our three-case strategy, and verifies lines~1
and~2 of Table~\ref{T:GVAWB}. Section~\ref{S:case1} considers Case 1
which includes $p=2$,
and Sections~\ref{S:case2} and~\ref{S:case3} consider Cases~2 and~3 when
$p$ is odd, and $G$ must be a nonabelian group of exponent~$p$.
Section~\ref{S:Ex} elucidates line~7 of Table~\ref{T:GVAWB} devotes substantial effort to constructing examples of
infinite and finite $k$-orbit groups for small $k$. The examples in line~6
of Table~\ref{T:GVAWB}
generalize to infinite 3- and 4-orbit groups in Lemma~\ref{L:GL3}.
Similarly, Example~\ref{E:Gpr} generalizes line~2 of Table~\ref{T:GVAWB} and
Lemmas~\ref{L1}, \ref{L:SpTr} generalize line~8. 
Finally, Section~\ref{S:4orbit} investigates the feasibility of classifying
4-orbit groups; as there are so few examples, a classification may be feasible.

\section{Preparation for a proof of Theorem~\texorpdfstring{\ref{T:main}}{}}
\label{S:prep}

In this section $G$ is a finite 3-orbit group. Let $N=\langle
G',\Phi(G)\rangle$, and let $A,B$ be the linear groups induced on
$G/N$ and $N$ by $\Aut(G)$ as in the Introduction.  Recall that
the \emph{solvable residual} $A^\infty$ of a finite group $A$, has
the property that $A/A^\infty$ is the largest solvable
factor group of $A$.
  We prove later (in Lemma~\ref{L:action})
  that $B^\infty\ne1$ implies $A^\infty\ne1$.
We split the proof of Theorem~\ref{T:main} into cases:
\[
  \textup{{\bf Case 1}.\;$A^\infty=1$,\ {\bf Case 2}.\;$B^\infty\ne1$ (so
    $A^\infty\ne1$), and\ {\bf Case 3}.\;$A^\infty\ne1$ and $B^\infty=1$.}
\]
In Case~1, $A$ is solvable. If $G$ is a $p$-group, then $\Aut(G)$ is
  solvable since the kernel
of the epimorphism $\Aut(G)\to A$ is the subgroup
$\textup{CAut}(G)\cong\Hom(V,W)\cong(\Cyc_p)^{mn}$ of central
automorphisms, see Lemma~\ref{L:action}.  In Theorem~\ref{T:Ainftyeq1}, $A$ is solvable or $p=2$.

The following lemmas can sometimes help to compute $\Aut(G)$
and $\omega(G)$. It uses the fact: If $M$ is characteristic in $G$
and $G>M>1$, then $\omega(G)\ge\omega(G/M)+\omega(M)-1$.

\begin{lemma}\label{L:C}  
  Suppose $G$ is a finite group and $G=M_1>\cdots>M_k=\{1\}$ where
  each $M_i$ is characteristic in $G$.
  Then $\omega(G)\ge 2-k+\sum_{i=1}^{k-1}\omega(M_i/M_{i+1})\ge k$ for $k\ge1$.
  If $\omega(G)=k$, then $G$ is solvable and has precisely $k$
  characteristic subgroups.
\end{lemma}

\begin{proof}
  The inequality $\omega(G)\ge 2-k+\sum_{i=1}^{k-1}\omega(M_i/M_{i+1})$ follows
  by induction on $k$ from $\omega(G)\ge\omega(G/M)+\omega(M)-1$.
  This show that $\omega(G)\ge k$ since
  $\omega(M_i/M_{i+1})\ge2$ for each $i$. Suppose that $\omega(G)=k$,
  i.e. equality holds. Then $\omega(M_i/M_{i+1})=2$ for each $i$, so each section
  $M_i/M_{i+1}$ is elementary abelian; whence $G$ is solvable. If $G$ has a
  characteristic subgroup different to the $k$ subgroups $M_1,\dots,M_k$,
  then it will refine the characteristic series $G=M_1>\cdots>M_k=\{1\}$, so
  $\omega(G)>k$, a contradiction.
\end{proof}



\begin{lemma}\label{L:A}
  If we know a lower bound $\omega(G)\ge\omega_0$ and a subgroup $A_0$ of
  $\Aut(G)$ with precisely $\omega_0$ orbits on $G$, then $\omega(G)\le\omega_0$
  and hence $\omega(G)=\omega_0$.
\end{lemma}

A finite nonabelian 3-orbit $p$-group has $1<G'\le\Phi(G)$ so
$N=G'=\Z(G)=\Phi(G)$. We now show that $A\le\GL(G/N)$ determines $B\le\GL(N)$,
and the $A$-orbits
on $V\cong G/N$ determine the $\Aut(G)$-orbits on~$G\setminus N$.

\begin{lemma}\label{L:action}
  For a nonabelian $3$-orbit $p$-group $G$, the preimage under
  the natural map $G\to G/\Z(G)\cong V$ of the nonzero $A$-orbits on $V$
  are the $\Aut(G)$-orbits on $G\setminus\Z(G)$.
  Further, the homomorphism $\phi\colon\Aut(G)\to\Aut(G/\Z(G))$ has image $A$ and
  kernel the central automorphisms $\textup{CAut}(G)\cong\textup{Hom}(V,W)$,
  and $A/K\cong B$ for some~$K\lhdeq A$.
\end{lemma}

\begin{proof}
  By definition, $\ker(\phi)$ equals $\textup{CAut}(G)$ namely the set of
  all $\alpha\in\Aut(G)$ that act trivially on $G/\Z(G)$.
  However, $G/\Z(G)\cong V$ and $\Z(G)\cong W$, and therefore
  $\textup{CAut}(G)\cong\textup{Hom}(V,W)$
  is elementary abelian of order~$|W|^{\dim(V)}=p^{mn}$.
  Furthermore, each $\alpha\in \textup{CAut}(G)$ acts trivially
  on $G'=\Z(G)\cong W$ since
  \[
    [g_1z_1,g_2z_2]^\alpha=[g_1^\alpha z_1^\alpha,g_2^\alpha z_2^\alpha]
    =[g_1^\alpha,g_2^\alpha]=[g_1,g_2]
    \qquad\textup{($g_1,g_2\in G, z_1,z_2\in\Z(G)$).}
  \]
  Hence the $A$-action on $V$ induces an action on $W$. The kernel of this
  action, namely $K\coloneq C_A(W)$, satisfies $\textup{CAut}(G)\le K$
  and $A/K\cong B$.
  The first sentence of the lemma is true since $g_1\in g_2\Z(G)$ implies $g_1^\alpha=g_2$
  for some $\alpha\in\textup{CAut}(G)$.
\end{proof}

The following straightforward lemma can be a useful tool for
computing $\Aut(G)$ where $G$ is a non-abelian 3-orbit $p$-group.
If we know that ${\mathcal A}_0\le\Aut(G)$ induces a subgroup $A_0$ of
$A=\Aut(G)^V$, $V=G/\Phi(G)$, and we show that $\Aut(G)$ can not induce
a larger subgroup of $\GL(V)$. Then $A=A_0$,
and hence $\Aut(G)=\CAut(G).{\mathcal A}_0$.

\begin{lemma}\label{L:Aut}
  Suppose we know a subgroup $A_0$ of $A\le\GL_m(p)$. If $A_0=\GL_m(p)$
  or $A_0$ is maximal (proper) subgroup of $\GL_m(p)$ and $A\ne\GL_m(p)$,
  then $A_0=A$.
\end{lemma}

Let $G$ be a finite 3-orbit group. If $G$ is abelian, then it is easy to
see that $\ord{G}\ne\{1,p\}$ or $\{1,p,r\}$ where $p\ne r$,
since $p,r\in\ord{G}$ implies $pr\in\ord{G}$. Hence
$\ord{G}=\{1,p,p^2\}$ so $G\cong (\Cyc_{p^2})^n\times(\Cyc_p)^k$ where
$\Cyc_{p^2}$ denotes a cyclic group of order~$p^2$ and $n\ge1$.
This implies $k=0$, otherwise $\omega(G)\ge4$ by Lemma~\ref{L:C} because of
the characteristic series
$G>\Omega_1(G)=\Cyc_p^{n+k}>\mho_1(G)=\Cyc_p^n>1$\label{P:line1}.
Thus line~1 of Table~\ref{T:GVAWB} is true.
Laffey and MacHale~\cite{LM} characterized the $3$-orbit groups $G$ with
$\ord{G}=\{1,p,r\}$ where $p\ne r$.
These groups are Frobenius groups of the form $W\rtimes V$ where $m=1$ and $r-1$
divides $n$, and they are related to projective geometry and
`uniform generation', see~\cite[Theorem~1.1]{G}.
This verifies line~2 of Table~\ref{T:GVAWB}. Certain $k$-orbit group
generalizations of line~2 (with $k\ge3$) are given in Example~\ref{E:Gpr}.
The `only' remaining case is when $G$ is a nonabelian $p$-group for
some prime $p$. This difficult case is not considered in~\cite{LM}, but is
covered in Sections~\ref{S:case1}--\ref{S:case3} and in~\cite{LZ2}.

\begin{hypothesis}\label{H:pgp}
  Let $G$ be a finite nonabelian 3-orbit $p$-group.
  Let $N=\Phi(G)$ and suppose that $G/N\cong V=(\F_p)^m$ and $N\cong W=(\F_p)^n$
  and $\Aut(G)$ induces subgroups $A\le\GL_m(p)$ and $B\le\GL_n(p)$ which act
  naturally and transitively on $V\setminus\{0\}$ and $W\setminus\{0\}$,
  respectively. Finally, let $K\lhdeq A$ be such that $A/K\cong B$
  as per Lemma~\ref{L:action}.
\end{hypothesis}

We assume Hypothesis~\ref{H:pgp} in Sections~\ref{S:case1}--\ref{S:case3}.
Thus $r=p$ and $N=G'=\Z(G)=\Phi(G)$ satisfies $1<N<G$.
Either $\exp(G)=p>2$, $\ord{G}=\{1,p\}$, and $\Aut(G)$ has two orbits on
elements of order $p$, or $\exp(G)=p^2$ and $\ord{G}=\{1,p,p^2\}$.
Hering's theorem~\cite{cH}
classifies the transitive linear subgroups $A\le\GL(V)$ and $B\le\GL(W)$;
our version is
Theorem~\ref{T:HeringNew}. The constraint $B\cong A/K$ 
(see Lemma~\ref{L:action}) further restricts the possibilities for $B$.
Our strategy is to compute possibilities for
the pair $(A,B)$ and then hopefully use the pair $(V,W)$ of modules to
reconstruct a unique 3-orbit group~$G$. Lemma~\ref{L:CentOdd}
describes how (and why) this is possible when $p>2$.

\section{Case 1 of Theorem~\texorpdfstring{\ref{T:main}}{}
  when \texorpdfstring{$A^\infty=1$}{} }
\label{S:case1}

  In this section we assume $p=2$ or $A^\infty=1$. We prove that $p=2$ and
  determine $G,A,B$.
The only 3-orbit group of order 8 is the quaternion group $Q_8$
and $\Aut(Q_8)\cong S_4$ is solvable. Our classification
is accelerated by appealing to~\cite[Proposition~3.1]{BG} which proves that
$\Aut(G)$ is solvable for any 3-orbit 2-groups $G$.

The possibilities for $G$
when $p=2$ were determined first in~\cite[Theorem~1.1]{BG} and then
in~\cite[Corollary~1.3]{LZ1}. The values of
$A$ and $B$ in lines of 3, 4, 5 of Table~\ref{T:GVAWB} of
Theorem~\ref{T:Ainftyeq1}
are proved in Lemma~\ref{L:abc}. Certainly $A$ follows from~\cite[Table~1]{LZ2}.
Cases~2,~3 relate to lines of~6,~7 of Table~\ref{T:GVAWB} and
are considered in Sections~\ref{S:case2},~\ref{S:case3}.

\begin{remark}\label{R:mgen}
  Let $G$ be a 3-orbit group.
  As $[g_1z_1,g_2z_2]=[g_1,g_2]$ for $z_1,z_2\in\Z(G)$, commutation gives rise
  to a bilinear map $V\times V\to W\colon (g_1\Z(G),g_2\Z(G))\mapsto[g_1,g_2]$.
  As $\omega(G)=3$, this map is surjective, so $|V|^2\ge|W|$
  or $2m\ge n$. We now prove the stronger bound $m\ge n$.
  If $p=2$, then the squaring map $Q\colon V\to W$ is surjective. Therefore
  $2^m=|V|\ge|W|=2^n$, and so $m\ge n$. Suppose now that $p>2$. If $n=1$, then
  $m\ge n$ holds, so suppose that $n\ge2$.
  Since $B$ is transitive on $W\setminus\{0\}$ we have
  \[
  p^n-1=|W\setminus\{0\}|\textup{ divides } |B|\textup{ divides }
  |A|\textup{ divides }|\GL_m(p)|.
  \]
  Whence $p^n-1$ divides $\prod_{i=1}^m(p^i-1)$. Since $n\ge2$ and $p>2$,
  Zsigmondy's theorem implies there exists a primitive prime divisor~$r$
  of $p^n-1$. As $r$ has order $n$ modulo~$p$ and $r$ divides
  $\prod_{i=1}^m(p^i-1)$, this shows that $m\ge n$, as claimed.\qed
\end{remark}

The history of, and properties of, Suzuki 2-groups
is summarized in~\cite[VIII.7]{HB}.

\begin{definition}\label{D:ABP}
  (a) Let $q=2^n$ and fix $\theta\in\Aut(\F_q)$ where $|\theta|>1$ is odd
  and $n\ne2^\ell$.
  The set $A(n,\theta)=\F_q\times\F_q$ with multiplication rule
  $(\lambda_1,\zeta_1)(\lambda_2,\zeta_2)=
  (\lambda_1+\lambda_2,\zeta_1+\zeta_2+\lambda_1^\theta\lambda_2)$
  defines a group of order~$q^2=2^{2n}$ called a
  \emph{Suzuki $2$-group of type A}.

  (b) Let $q=2^n$ where $n\ge1$ and\footnote{Higman~\cite{Hig} had $n>1$, but we will
  allow $n=1$ and $q=2$ to include the quaternion group $Q_8$.} fix $\eps\in\F_{q^2}^\times$ of order $q+1$.
  The set $B(n)=\F_{q^2}\times\F_q$ with multiplication rule
  $(\lambda_1,\zeta_1)(\lambda_2,\zeta_2)=
  (\lambda_1+\lambda_2,\zeta_1+\zeta_2
    +\lambda_1\lambda_2^q\eps+(\lambda_1\lambda_2^q\eps)^q)$
  defines a group of order $q^3$ whose isomorphism type is
  independent of $\eps$, see~\cite[Theorem~(v)]{Dornhoff}.

  (c) Let $q=2^3$ and fix $\eps\in\F_{q^2}^\times$ of order $q^2-1$.
  The set $P=\F_{q^2}\times\F_q$ with multiplication rule
  $(\lambda_1,\zeta_1)(\lambda_2,\zeta_2)=
  (\lambda_1+\lambda_2,\zeta_1+\zeta_2
    +\lambda_1\lambda_2^2\eps+(\lambda_1\lambda_2^2\eps)^q)$
  defines a group of order $q^3=2^9$ with isomorphism type
  independent of $\eps$, see~\cite[p.\,705]{Dornhoff}.
\end{definition}

\begin{theorem}\label{T:Ainftyeq1}
  Let $G$ be a finite nonabelian $3$-orbit $p$-group with $|\Phi(G)|=p^n$ and
  let $V,A,W,B$ be as in Hypothesis~$\ref{H:pgp}$. If $p=2$ or $A^\infty=1$, then
  $G,V,A,W,B$ are as in lines $3, 4, 5$ of Table~$\ref{T:GVAWB}$. In particular,
  $p=2$ and $\Aut(G)$ is solvable.
\end{theorem}

\begin{proof}
  If $p=2$, then $\Aut(G)$ is solvable by~\cite[Proposition~3.1]{BG}.
  Assume now that $A^\infty=1$, so that $A$ is solvable. Since
  $\Phi(G)=\Z(G)$, the kernel of the homomorphism $\Aut(G)\to A$ is
  abelian by Lemma~\ref{L:action}, so $\Aut(G)$ is solvable. Thus $G$
  is listed in~\cite[Theorem]{Dornhoff}.  The nonabelian $p$-groups in
  this list are in cases~(iv)--(viii) of~\cite[Theorem]{Dornhoff}.
  The groups in cases (vi), (vii) are excluded because
  $\Aut(G)^\infty\ne1$, see Lemma~\ref{L:SpTr}. The remaining groups
  in cases (iv), (v), (viii) are those in
  Definition~$\ref{D:ABP}${\rm(a,b,c)} including the quaternion group
  $Q_8$ which is $B(1)$. Hence $p=2$.  The values of $A, B$ (and $V,
  W$) for the remaining cases follow from Lemma~\ref{L:abc}(a,b,c).
\end{proof}

\section{Case 2 of Theorem~\texorpdfstring{\ref{T:main}}{}
  when \texorpdfstring{$B^\infty\ne1$}{} }
\label{S:case2}

In this section $G$ is a nonabelian finite 3-orbit $p$-group where $p$ is odd,
and we assume that $B$ is not solvable. Hence the solvable
residual $B^\infty$ of $B$ is nontrivial. Our classification of possible $G$
relies on Hering's theorem~\cite{cH},
which is proved in~\cite[Appendix~A]{Lie87a}, and classifies the subgroups
$A\le\GL_m(p)$ which act transitively on the nonzero vectors of the natural
module $V=\F_p^{\,m}$. The following more concise version of Hering's theorem
constrains $A^\infty$ instead of $A$. It follows easily
from~\cite[Appendix~A]{Lie87a}, and we let the reader check the details.
(We used {\sc Magma}~\cite{Magma} to find $H$ in part~(a).)

\begin{theorem}[Hering~\cite{cH}]\label{T:HeringNew}
  Let $A\le\GL_m(p)$ act transitively on the nonzero
  vectors of the natural module $\F_p^{\;m}$. Then the solvable residual
  $A^\infty$ of $A$ lies in the list:
  \begin{enumerate}[{\rm (a)}]
  \item $A^\infty=1$ if $A\le\GammaL_1(p^m)$ or if $m=2$, $Q_8\lhd A\le (Q_8.S_3)\circ\Cyc_{p-1}$ and $p\in\{5,7,11,23\}$ or if $(m,p)=(4,3)$, $A=(D_8\circ Q_8).H$ where $H\in\{\Cyc_5,D_{10},F_{20}\}$,
  \item $A^\infty=\begin{cases}
      \SL_{m/b}(p^b)&\textup{if $2\le m/b\le m$ and $(m/b,p^b)\ne(2,3)$,}\\
      \Sp_{m/b}(p^b)&\textup{if $m/b\ge4$ is even,}\\
      \textup{G}_2(2^b)'&\textup{if $(m,p)=(6b,2)$,}\end{cases}$
  \item $A^\infty=\SL_2(5)$ where $(m,p)\in\{(4,3),(2,11),(2,19),(2,29),(2,59)\}$,  \item $A^\infty=A=\begin{cases}
      \textup{$A_6$ or $A_7$}&\textup{if $(m,p)=(4,2)$,}\\
      \SL_2(13)&\textup{if $(m,p)=(6,3)$.}\end{cases}$
  \item $A^\infty=2^{1+4}_{-}.A_5$ where $(m,p)=(4,3)$
    and $2^{1+4}_{-}\cong D_8\circ Q_8$ is extraspecial.
  \end{enumerate}
  Furthermore, $A^\infty$ is transitive on
    the nonzero vectors $\F_p^{\;m}\setminus\{0\}$ in cases \textup{(b,c,d,e)}.
\end{theorem}

The Fitting subgroup, $\Fit(A)$, of $A$ is the largest normal nilpotent subgroup of $A$.



\begin{corollary}\label{C:Hering}
  If $A\le\GL_m(p)$ is not solvable and acts transitively on the nonzero
  vectors of the natural module $\F_p^{\;m}$, then $A^\infty$
  is also transitive on $\F_p^{\;m}\setminus\{0\}$ and
  $A^\infty/\Fit(A^\infty)$ is a nonabelian simple group.
  Further, if $B$ is not solvable and acts transitively on
  $\F_p^{\;n}\setminus\{0\}$, then 
  $B^\infty/\Fit(B^\infty)\cong A^\infty/\Fit(A^\infty)$ is nonabelian and simple.
\end{corollary}

\begin{proof}
    Now $A^\infty\ne1$ acts transitively on $\F_p^{\;m}\setminus\{0\}$ by
    Theorem~\ref{T:HeringNew}, and $A^\infty/\Fit(A^\infty)$
  is a nonabelian simple group.
  By Lemma~\ref{L:action}, $B$ is a quotient of $A$, so $B^\infty$ is
  a quotient of $A^\infty$. If $B$ is not solvable and acts transitively
  on $\F_p^{\;n}\setminus\{0\}$, then $B^\infty\ne1$ so $A^\infty\ne1$ and
  $B^\infty$ is also transitive on $\F_p^{\;n}\setminus\{0\}$
  by Theorem~\ref{T:HeringNew}.
  In addition, Theorem~\ref{T:HeringNew}(b,c,d,e) implies that
    $B^\infty/\Fit(B^\infty)\cong A^\infty/\Fit(A^\infty)$ is nonabelian and simple.
\end{proof}

\begin{lemma}\label{L:CentOdd}
  Suppose $G,p,V,A,W,B$ are as in Hypothesis~$\ref{H:pgp}$ where $p$ is odd.
  Then
  \begin{enumerate}[{\rm (a)}]
  \item $G$ has exponent $p$ and hence is a factor group of 
    group $V\ltimes\Lambda^2V$ and $n\le\binom{m}{2}$.
  \item The center $\Z(A)$ is cyclic and
    $|\Z(A)|$ divides $p^{e_A}-1$ where $e_A\mid m$. 
  \item If $B^\infty\ne1$, then $|\Z(A^\infty)|$ is odd.
  \end{enumerate}
\end{lemma}

\begin{proof}
  (a)~By Hypothesis~$\ref{H:pgp}$, $G$ is an $m$-generated nonabelian
  3-orbit $p$-group. If $G$ has exponent $p^2$, then $\Aut(G)$ is transitive
  on the elements of order $p$. Since $p>2$, $G$ is abelian by a deep result of
  Shult~\cite[Corollary~3]{eS}. This contradiction shows that $G$ has
  exponent $p$. Thus $G$ is a factor group of the (universal) $m$-generated
  exponent-$p$ class 2 group, the elements of which may be viewed as
  ordered pairs $(v,w)\in V\times\Lambda^2V$ with multiplication rule
  $(v_1,w_1)(v_2,w_2)=(v_1+v_2,w_1+w_2+v_1\wedge v_2)$. Thus $n\le\binom{m}{2}$.
  
  (b)~As $A$  acts irreducibly on $V$, results of Schur and Wedderburn
  imply that the ring ${\rm End}_{\F_pA}(V)$ of $\F_pA$-endomorphisms
  is a finite field, say $\F_q$ where $q=p^{e_A}$ depends on~$A$. Further,
  $\Z(A)\le\F_q^\times$ is cyclic and $|\F_q^\times|=p^{e_A}-1$ where $e_A\mid m$.
  
  (c)~Since $B^\infty\ne1$ and $A/K\cong B$ (Lemma~\ref{L:action}),
  we have $A^\infty\ne1$.
  Hence $A^\infty/\Fit(A^\infty)\cong B^\infty/\Fit(B^\infty)$
  is the unique nonabelian simple composition factor of $A$ and $B$ by
  Corollary~\ref{C:Hering}.
    By Theorem~\ref{T:HeringNew}, $A^\infty$ acts transitively on
    $\F_p^{\,m}\setminus\{0\}$, so by part~(b), $\Z(A^\infty)$ is cyclic
    and $|\Z(A^\infty)|$ divides $p^m-1$.
  If $|\Z(A^\infty)|$ is even, then $-1\in\Z(A^\infty)$
  by part~(b). However $-1$ acts trivially on
  $\Lambda^2V$ and hence on $W=\Lambda^2V/U$. Therefore $B^\infty\ne A^\infty$,
  contradicting Theorem~\ref{T:HeringNew}. Thus $|\Z(A^\infty)|$~is~odd.
\end{proof}

We remark that constructing 3-orbit $p$-groups of (odd)
exponent $p$ is the same as finding
maximal $A$-submodules of the exterior square of an $A$-module.

\begin{remark}\label{R:maximal}
  In Lemma~\ref{L:CentOdd}(a), $A$ acts on $V$ and hence on $\Lambda^2V$.
  As $G$ is a 3-orbit group, $A$ acts irreducibly on the quotient $\Lambda^2V/U$
  of $\Lambda^2V$, so $U$
  is a maximal $A$-submodule; and $B$ acts faithfully on $\Lambda^2V/U$.
  The group
  ${\mathcal G}\coloneq V\times\Lambda^2V$ has center $\{0\}\times\Lambda^2V$,
  hence $U\lhd {\mathcal G}$ and $G\cong {\mathcal G}/U$. Thus
  $G\cong V\times(\Lambda^2V/U)$ has multiplication rule
  $(v_1,w_1+U)(v_2,w_2+U)=(v_1+v_2,w_1+w_2+v_1\wedge v_2+U)$.
  If $\alpha\in A$ and $\alpha K\in B$, then $\alpha$ acts as
  $(v_1,w_1+U)^\alpha=(v_1^\alpha,(w_1+U)^{\alpha K})$.
  Thus  a 3-orbit group $G$ gives rise to a maximal $A$-submodule $U$
  of $\Lambda^2V$. Conversely, $A$ may not be transitive on the nonzero vectors
  of $\Lambda^2V/U$ where $U$ is a maximal $A$-submodule
  of $\Lambda^2V$. Interestingly, this is not the case,
  see Remark~\ref{R:sub}.\qed
\end{remark}

The following fact, follows from~\cite[Line 2, Table~3]{Lie87a} when $d\ge3$.
\begin{equation}\label{E:SL3}
  \textup{If $\V=\F_q^{\,d}$ is the natural module for $\SL_d(q)$, then
    $\Lambda^2\V$ is irreducible.}
\end{equation}
Now $\dim(\Lambda^2\V)=\binom{d}{2}$, so $\Lambda^2\V$ is a 1-dimensional
trivial module if $d=2$. If $d=3$, then $\Lambda^2\V$ is isomorphic to the dual
module $\V^*$ of $\V$. If $q=p^b$, then the $b$ Galois conjugate modules
$\V^\theta$, $\theta\in\Gal(\F_q/\F_p)$, all give rise to the same
irreducible $db$-dimensional $\F_p\SL_d(q)$-module by~\cite[VII.1.16]{HB}.
This process of changing from $\V$ to $V=\V\restrict F=\F_p^{\,bd}$ is
sometimes called
`blowing up the dimension' or `restricting to a subfield'.

\begin{theorem}\label{T:Binftyne1}
  Let $G$, $p$, $m$, $V$, $A$, $n$, $W$, $B$ be as in Hypothesis~$\ref{H:pgp}$. 
  If $B^\infty\ne1$, then $p$ is an odd prime, $m=n$ is divisible by $3$ and
  $A^\infty\cong B^\infty\cong\SL_3(p^{n/3})$ where $A^\infty$ acts as an
  $(\F_{p^{n/3}})^3$-module, and $B^\infty$ acts as its dual. Furthermore,
  $G$ is isomorphic to the class~$2$ factor group in Lemma~$\ref{L:GL3}$ with
  $F=F_0=\F_{p^{n/3}}$, as on line~$6$ of Table~$\ref{T:GVAWB}$.
\end{theorem}

\begin{proof}
  By Corollary~\ref{C:Hering},
  $B^\infty/\Fit(B^\infty)\cong A^\infty/\Fit(A^\infty)$ is a nonabelian
  simple group. Recall that $G$ is nonabelian by Hypothesis~$\ref{H:pgp}$.
  Theorem~\ref{T:Ainftyeq1} implies that $p\ne2$ otherwise $A^\infty=B^\infty=1$.
  Thus $n\le\binom{m}{2}$ by Lemma~\ref{L:CentOdd}(a). This shows
  that $m\ne2$, otherwise $n=1$ and $B^\infty=1$ as $\GL_1(p)$ is cyclic.
  Hence $m\ge3$ and $p>2$. We can rule out case~(a) of
  Theorem~\ref{T:HeringNew}, and $\textup{G}_2(2^b)'$ in case~(b) as $p\ne2$.

  As $|\Z(A^\infty)|$ is odd (Lemma~\ref{L:CentOdd}(c)), 
  $A^\infty\not\in\{\Sp_{m/b}(p^b),\SL_2(5),\SL_2(13),2^{1+4}_{-}.A_5\}$ in
  Theorem~\ref{T:HeringNew}. This rules out cases (c), (d) and (e)
  of Theorem~\ref{T:HeringNew}.
  Thus  $A^\infty=\SL_{m/b}(p^b)$ and $m/b$ is odd
  because $|\Z(\SL_{m/b}(p^b))|=\gcd(m/b,p^b-1)$ is odd. Hence
  $A^\infty\cong B^\infty$ by Theorem~\ref{T:HeringNew}, so $m=n$. Let $V$ be
  the natural $m$-dimensional $A^\infty$-module over $\F_p$, and let $\V$ be
  the $d$-dimensional $A^\infty$-module over $\F_q$ where $q=p^b$.
  The exterior square $\Lambda^2\V$ has dimension $\binom{d}{2}$ over the field
  $\F_q\cong\textup{End}_{\F_pA^\infty}\V$ of endomorphisms.
  However, $\V$ is irreducible and so too is $\Lambda^2\V$ by~\eqref{E:SL3}.
  Hence by~\cite[VII.1.16(e)]{HB} (see also Remark~\ref{R:sub}), we have
  $\binom{m/b}{2}=n/b=m/b$, so $m/b=3$ and $m=3b=n$.

  In summary, $A^\infty\cong\SL_{3}(p^b)$ acts faithfully on the 3-space $\V$ and
  $B^\infty\cong\SL_{3}(p^b)$ acts faithfully and irreducibly
  on the 3-space $\Lambda^2\V$.
  Adapting the argument in Remark~\ref{R:maximal}, $G$ is isomorphic
  to the group $\V\times\Lambda^2\V$.
  Alternatively, we may identify $G$ with the set $\F_q^{\,3}\times\F_q^{\,3}$
  with multiplication given in Lemma~\ref{L:GL3} (by ignoring the
  third coordinate). Setting $F=F_0=\F_q$ in Lemma~\ref{L:GL3} shows
  that $A$ contains
  $\GammaL_3(\F_q)=\GL_3(\F_q)\rtimes\Aut(\F_q)$. If $q=p^b$, then
  $\GammaL_3(\F_q)$ is a maximal proper subgroup of $\GL_{3b}(\F_p)$ if $b>1$,
  so Lemma~\ref{L:Aut} implies that $A=\GammaL_3(\F_q)$ for $b\ge1$.
  A similar argument shows that
  $B=\GammaL^{+}_3(\F_q)=\{g\in\GL_3(\F_q)\mid\det(g)\in(\F_q^\times)^2\}\rtimes\Aut(\F_q)$
  by Remark~\ref{R:det} and Lemma~\ref{L:GL3}.
  This verifies line~6 of Table~\ref{T:GVAWB}.
\end{proof}

Nonabelian $p$-groups with precisely 3 characteristic subgroups
were called UCS groups (Unique Characteristic Subgroup) by Taunt
and were studied in~\cite{GPS}.
A 3-orbit $p$-group $G$ is a UCS group which is a special group
as $\Z(G)=G'=\Phi(G)$.
The structure of a special group $G$ is strongly influenced by the
two $\Aut(G)$-modules $V=G/\Z(G)$ and $W=\Z(G)$, see~\cite[Theorem~6]{GPS}.
The group $G_4$ in~\cite[Theorem~8]{GPS} is an exponent-$p^2$ cousin
(with $B=\textup{SO}_2(p)$)
of the exponent-$p$ groups in Theorem~\ref{T:Binftyne1}.

\section{Case 3 of Theorem~\texorpdfstring{\ref{T:main}}{}:
  \texorpdfstring{$A^\infty\ne1$}{} and \texorpdfstring{$B^\infty=1$}{} }
\label{S:case3}

In this section $G$ is a 3-orbit $p$-group where $p$ is odd, and
$A^\infty\ne1$ and $B^\infty=1$ hold.
A prototypical example is an extraspecial $p$-group as hinted by the following
lemma. It is worth examining this case before considering more general
examples.

The 3-orbit group $Q_8$ was considered in Section~\ref{S:case1}.
We show that extraspecial 2-groups are a rich source
of 4-orbit groups. We focus primarily on the case $p>2$.

\begin{lemma}\label{L1}
  A finite extraspecial $p$-group $G$ is a $3$-orbit group precisely when
  $G\cong Q_8$ or $G\cong p_{+}^{1+m}$ has odd exponent~$p$. In these
  cases $\Aut(G)$ induces on $\Z(G)$ the cyclic subgroup $B\cong\Cyc_{p-1}$.
  An extraspecial $2$-group $G$ with $G\not\cong Q_8$ is a $4$-orbit group.
\end{lemma}

\begin{proof}
  A finite extraspecial $p$-group $G$ satisfies $G'=\Z(G)=\Phi(G)\cong\Cyc_p$.
  Suppose that $G$ is a 3-orbit group and the notation in Hypothesis~\ref{H:pgp}
  holds. If  $G\setminus\Z(G)$  contains elements of orders $p$ and $p^2$, then
  $\Aut(G)$ has at least 4 orbits on $G$. This is the case if $G$ has
  odd exponent $p^2$, or $p=2$ and $G\not\cong Q_8$. However, if $G\cong Q_8$,
  then $\Aut(G)$ induces $\GL_2(2)$ on $V=G/\Z(G)\cong(\Cyc_2)^2$, and acts
  trivially on $\Z(G)\cong\Cyc_2$.
  Hence $Q_8$ is a $3$-orbit group. It follows from~\cite{W} that
  an extraspecial $p$-group of odd exponent~$p$ is a 3-orbit group.
  This is also proved in
  Lemma~\ref{L:SpTr} which also applies to infinite $3$-orbit groups.
  In our case $B\cong\GL_1(p)\cong\Cyc_{p-1}$ holds by Lemma~\ref{L:SpTr}.

  Suppose now that $G_\eps$ is the extraspecial 2-group $2_\eps^{1+m}$
  of order $2^{m+1}$ and type $\eps\in\{-,+\}$ where $m$ is even.
  In $G_\eps$, squaring induces a (well-defined) quadratic form $Q_\eps$
  on the vector space $V_\eps=G_\eps/\Z(G_\eps)\cong\F_2^{\,m}$. The preimage
  in $G_\eps$ of singular vectors in $V_\eps$ are the noncentral involutions of
  $G_\eps$, and the preimage of nonsingular vectors in $V_\eps$ are the
  elements of order 4 in $G_\eps$. It is well known
  that the outer automorphism group $\textup{Out}(G_\eps)$ is isomorphic
  to the full orthogonal group ${\rm O}(Q_\eps)\cong{\rm O}_m^\eps(2)$,
  see~\cite{G95} or \cite[\S2.2.6]{BHRD}. 
  For even $m\ge2$ and $\eps\in\{-,+\}$ the space $V_\eps$ has nonsingular
  vectors,
  and it has singular vectors except when $(\eps,m)=(-,2)$. By Witt's theorem
  ${\rm O}_m^\eps(2)$ is transitive on the (possibly empty) set of singular
  vectors and the set of nonsingular vectors. Clearly $\{1\}$ and
  $\Z(G_\eps)\setminus\{1\}$ are $\Aut(G_\eps)$-orbits, so that $G_\eps$ is
  a 3-orbit group if $(\eps,m)=(-,2)$, and a 4-orbit group otherwise.
  The elements of order 4 form one $\Aut(G_\eps)$-orbit, and the involutions
  form two orbits if $(\eps,m)\ne(-,2)$ (the central involution is fixed).
\end{proof}

The following fact from representation theory will guide our proof
of Theorem~\ref{T:Binftyeq1}.

\begin{remark}\label{R:Sn}
  The symmetric group $S_n$ acts on $\F_q^{\,n}$ by permuting the elements of
  a basis $\{v_1,\dots,v_n\}$. Further, the augmentation map
  $\phi\colon V\to\F_q\colon\sum_{i=1}^n \lambda_iv_i\mapsto\sum_{i=1}^n \lambda_i$
  is an $S_n$-epimorphism, and $S_n$ fixes the submodules
  \[
  W=\ker(\phi)=\left\{\sum_{i=1}^n \lambda_iv_i\mid \sum_{i=1}^n\lambda_i=0\right\}
  \quad\textup{and}\quad D=\left\langle\sum_{i=1}^nv_i\right\rangle.
  \]
  The equation $\sum_{i=1}^ni(v_i-v_{i+1})=(\sum_{i=1}^{n-1}v_i)-(n-1)v_n$
  shows that $D\subseteq W$
  if $p\mid n$ where $p=\Char(\F_q)$, and $D\cap W=\{0\}$ otherwise. We call
  $W/(D\cap W)$ the fully deleted permutation module. It can be
  written over $\F_p$, and it is absolutely irreducible.\qed
\end{remark}

\begin{remark}\label{R:Spdq}
  Let $\V=\F_q^{\,2\ell}$ be a symplectic space preserving the
  non-degenerate alternating bilinear form $f\colon\V\times\V\to\F_q$. Let
  $\V$ be the natural module for  $\Sp(\V)\cong\Sp_{2\ell}(q)$ where
  $q=p^b$, with basis $e_1,\dots,e_{2\ell}$. The map
  $\phi\colon\Lambda^2\V\to\F_q$ with
  $\sum\lambda_{ij}e_i\wedge e_j\mapsto\sum\lambda_{ij}f(e_i,e_j)$ is
  a $\Sp(\V)$-module epimorphism. Hence $\mathcal{W}\coloneq\ker(\phi)$ is an
  $\Sp(\V)$-invariant hyperplane.
  Let $\langle e_1,e_2\rangle,\dots,\langle e_{2\ell-1},e_{2\ell}\rangle$ be
  pairwise orthogonal hyperbolic planes. Then
  $\phi(\sum_{i<j}\lambda_{ij}e_i\wedge e_j)=\sum_{i=1}^\ell\lambda_{2i-1,2i}$ and
  the stabilizer $\Sp_2(q)\wr S_\ell$ of the decomposition
  $\langle e_1,e_2\rangle\oplus\cdots\oplus\langle e_{2\ell-1},e_{2\ell}\rangle$ is
  a maximal subgroup of $\Sp_{2\ell}(q)$ which preserves the submodules
  $\mathcal{W}=\ker(\phi)$ and $\mathcal{D}=\langle\sum_{i=1}^\ell e_{2i-1}\wedge e_{2i}\rangle$ in
  Figure~\ref{F:Spdq} by Remark~\ref{R:Sn}. A symplectic transvection not
  in $\Sp_2(q)\wr S_\ell$ also preserves these submodules. Hence $\mathcal{D}$ and $\mathcal{W}$ are
  invariant under all of $\Sp_{2\ell}(q)$.
\begin{figure}[!ht]
\begin{center}
\begin{tikzpicture}[scale=0.5]
\draw [-] (0,0) -- (0,4);
\draw [fill=black] (0,0) circle [radius=1mm];
\draw [fill=black] ((0,1) circle [radius=1mm];
\draw [fill=black] ((0,3) circle [radius=1mm];
\draw [fill=black] ((0,4) circle [radius=1mm];
\node [left,scale=1] at (-0.2,4.1) {$\Lambda^2\V$};
\node [left,scale=1] at (-0.2,2.9) {$\mathcal{W}$};
\node [left,scale=1] at (-0.2,1.1) {$\mathcal{D}$};
\node [left,scale=1] at (-0.2,0) {$\{0\}$};
\node [right,scale=1] at (0,0.5) {$1$};
\node [right,scale=1] at (0,2) {$2\ell^2-\ell-2$};
\node [right,scale=1] at (0,3.5) {$1$};
\node [below,scale=1] at (0,-0.5) {$p$ divides $\ell$};
\end{tikzpicture}
\hskip10mm
\begin{tikzpicture}[scale=0.44]
\draw [-] (0,0) -- (0,1) -- (2,4) -- (2,3) -- (0,0);
\draw [fill=black] (0,0) circle [radius=1mm];
\draw [fill=black] (0,1) circle [radius=1mm];
\draw [fill=black] (2,3) circle [radius=1mm];
\draw [fill=black] (2,4) circle [radius=1mm];
\node [right,scale=1] at (2.2,4.1) {$\Lambda^2\V$};
\node [right,scale=1] at (2.2,2.9) {$\mathcal{W}$};
\node [left,scale=1] at (0,1.4) {$\mathcal{D}$};
\node [right,scale=1] at (0.2,0) {$\{0\}$};
\node [left,scale=1] at (0,0.5) {$1$};
\node [right,scale=1] at (1.1,1.5) {$2\ell^2-\ell-1$};
\node [below,scale=1] at (0,-0.5) {\qquad\quad$p$ coprime to $\ell$};
\end{tikzpicture}
\end{center}
\vskip-6mm
\caption{The $A$-submodules of $\Lambda^2\V$ where $\V=\F_{p^b}^{\,2\ell}$ and
$\Sp(\V)\le A\le\Gamma\Sp(\V)$}\label{F:Spdq}
\end{figure}

  We show that $\mathcal{D}$ and $\mathcal{W}$ are $A$-invariant for $A$
  satisfying $\Sp(\V)\le A\le\Gamma\Sp(\V)$. The notation $\CSp(\V)$ and
  $\Gamma\Sp(\V)$ is described in Remark~\ref{R:Sp}.
  First, $\CSp_{2\ell}(q)=\langle g_\mu,\Sp_{2\ell}(q)\rangle$ where
  $\F_q^\times=\langle\mu\rangle$ and $g_\mu$ satisfies
  $e_{2i-1}g_\mu=\mu e_{2i-1}$ and $e_{2i}g_\mu=e_{2i}$ for $i\le\ell$. 
  Also $\phi(ug_\mu)=\mu\phi(u)$ for $u\in\Lambda^2\V$, so $\CSp_{2\ell}(q)$ fixes
  $\mathcal{W}=\ker(\phi)$ and $\mathcal{D}$. Second, if $g_\theta$ satisfies
  $(\sum_{i=1}^{2\ell}\lambda_ie_i)g_\theta=\sum_{i=1}^{2\ell}\lambda_i^\theta e_i$ for
  $\theta\in\Aut(\F_{p^b})$, then $\phi(ug_\theta)=\phi(u)^\theta$
  for $u\in\Lambda^2\V$. Hence
  $\Gamma\Sp_{2\ell}(p^b)$ fixes $\mathcal{W}$ and $\mathcal{D}$ and
  induces $\GammaL_1(p^b)$ on $\mathcal{D}$.
  Finally, the only $\Sp(\V)$-submodules of $\Lambda^2\V$ are
  $\{0\}$, $\mathcal{W}$, $\mathcal{D}$, $\Lambda^2\V$
  by~\cite[Table~5]{Lie87a}, see also~\cite[Hauptsatz~1]{Hiss} when~$p=2$.
  In summary, we have justified Figure~\ref{F:Spdq}.\qed
\end{remark}

\begin{remark}\label{R:AS}
  Let $V=F^m$ be an $m$-dimensional vector space over a field $F$
  where $\Char(F)\ne2$. Then $T^2(V)=A^2(V)\oplus S^2(V)$ where 
  $T^2(V)=V\otimes V$, $A^2(V)$ and $S^2(V)$ are called the tensor, alternating,
  and symmetric squares of $V$, respectively. Set
  $A^2(V)\coloneq\{v_1\otimes v_2-v_2\otimes v_1\mid v_1,v_2\in V\}$ and
  $S^2(V)\coloneq\{v_1\otimes v_2+v_2\otimes v_1\mid v_1,v_2\in V\}$.
  The identity
  $v_1\otimes v_2=
  \frac12(v_1\otimes v_2-v_2\otimes v_1)+
  \frac12(v_1\otimes v_2+v_2\otimes v_1)$ for $v_1,v_2\in V$ implies that
  $T^2(V)=A^2(V)\oplus S^2(V)$ holds.
  The exterior square is isomorphic to the alternating square
  via $\Lambda^2(V)\to A^2(V)
  \colon v_1\wedge v_2\mapsto \frac12(v_1\otimes v_2-v_2\otimes v_1)$.
  The symmetric square is normally defined to be the quotient
  $T^2(V)/A^2(V)$ (in all characteristics), and similarly $\Lambda^2(V)$
  is defined to be $T^2(V)/S^2(V)$ in all characteristics. The isomorphism
  $T^2(V)/A^2(V)\to S^2(V)$ with $v_1\otimes v_2+A^2(V)\mapsto
  v_1\otimes v_2+v_2\otimes v_1$ holds for all $F$.\qed
\end{remark}

\begin{remark}\label{R:C}
  Let $p$ be an odd prime and set $q=p^b$, $E=\F_q$, $F=\F_p$, $\V=E^d$ where
  $d=2\ell$ is even and $\V$ is an $A$-module
  where $\Sp(\V)\le A\le\Gamma\Sp(\V)$.
  We view $V=\V\restrict F=F^{bd}$ as the $A$-module $\V=E^d$
  written over $F$ using the inclusions $A\le\Gamma\Sp_d(E)\le\GL_{bd}(F)$.
  The $bd$-dimensional $EA$-module $V\otimes_F E$ is
  $\bigoplus_{i=0}^{b-1}\V^{(i)}$ where $\V^{(i)}$ it the Galois conjugate of
  $\V$ by $\theta^i$ where $\Gal(E:F)=\langle\theta\rangle$
  by~\cite[VII.1.16(a)]{HB}. We shall prove the decomposition
  $A^2(V\otimes E)=\bigoplus_iA^2(\V^{(i)})\oplus\bigoplus_{i<j}\mathcal{A}_{ij}$
  where $0\le i<b$, $0\le i<j<b$, and $\mathcal{A}_{ij}$ is defined below.
  First, $T^2(V)\otimes E\cong T^2(V\otimes E)$ equals
  \[
  \left(\bigoplus_i\V^{(i)}\right)\otimes\left(\bigoplus_j\V^{(j)}\right)
  =\bigoplus_iT^2(\V^{(i)})\oplus
  \bigoplus_{i<j}\left(\V^{(i)}\otimes\V^{(j)}\oplus\V^{(j)}\otimes\V^{(i)}\right).
  \]
  However, $\V^{(i)}\otimes\V^{(j)}\oplus\V^{(j)}\otimes\V^{(i)}$ is a direct sum
  of submodules say $\mathcal{A}_{ij}\oplus\mathcal{S}_{ij}$ where 
  \begin{align*}
  \mathcal{A}_{ij}&=
  \left\{v_i\otimes v_j-v_j\otimes v_i\mid v_i\in\V^{(i)},v_j\in\V^{(j)}\right\},\\
  \mathcal{S}_{ij}&=
  \left\{v_i\otimes v_j+v_j\otimes v_i\mid v_i\in\V^{(i)},v_j\in\V^{(j)}\right\},
  \end{align*}
  $\mathcal{A}_{ij}\le A^2(V\otimes E)$ and
  $\mathcal{S}_{ij}\le S^2(V\otimes E)$ by Remark~\ref{R:AS}. Hence
  \[
  T^2(V\otimes E)=\bigoplus_iA^2(\V^{(i)})\oplus\bigoplus_iS^2(\V^{(i)})
  \oplus\bigoplus_{i<j}\mathcal{A}_{ij}\oplus\bigoplus_{i<j}\mathcal{S}_{ij}.
  \]
  Since
  $\mathcal{A}_{ij}\cong \V^{(i)}\otimes\V^{(j)}\cong\mathcal{S}_{ij}$ for $i<j$,
  the claimed decomposition follows: 
  \[
  A^2(V\otimes E)=\bigoplus_iA^2(\V^{(i)})\oplus\bigoplus_{i<j}\mathcal{A}_{ij}
  \cong\bigoplus_iA^2(\V^{(i)})\oplus\bigoplus_{i<j}\V^{(i)}\otimes\V^{(j)}.
  \]
  (The containment $\ge$ holds, and the dimension agree
  as $\binom{bd}{2}=b\binom{d}{2}+\binom{b}{2}d^2$.)\qed
\end{remark}

Recall that $A^\infty\ne1$ and $B^\infty=1$ hold in this section.

\begin{remark}\label{R:sub}
Assume, as in Remark~\ref{R:C}, that $p>2$ is prime, $q=p^b$,
$E=\F_q$, $F=\F_p$, $\V=E^d$, $d=2\ell$ is even, $V=\V\restrict F=F^{bd}$
where $\V$ is an $A$-module and $\Sp(\V)\le A\le\Gamma\Sp(\V)$.
By Remark~\ref{R:Spdq}, $E$ is a $b$-dimensional $F\kern-2pt A$-module which
is irreducible if $\CSp(\V)\le A$ and is trivial if $A=\Sp(\V)$. 
Let $U$ be a maximal $F\kern-2pt A$-submodule of
$A^2(V\otimes E)\restrict F$.  This remark proves that the quotient
$F\kern-2pt A$-module $(A^2(V\otimes E)\restrict F)/U$ is isomorphic
to a subfield of $E$ containing $F$ (really a quotient $F\kern-2pt A$-module).
As $B^\infty=1$, we see that
$A^\infty=\Sp(\V)$ acts trivially on this quotient. In the next paragraph, we
consider $EA$-submodules rather than $F\kern-2pt A$-submodules.

  Let $\mathcal{U}$ be 
  a maximal $EA$-submodule of $A^2(V\otimes E)$ such that $A^\infty=\Sp(\V)$ acts
  trivially on $A^2(V\otimes E)/\mathcal{U}$. 
  Remark~\ref{R:C} shows that $A^2(V\otimes E)=\mathcal{X}\oplus \mathcal{Y}$
  where $\mathcal{X}=\bigoplus_iA^2(\V^{(i)})$ and
  $\mathcal{Y}=\bigoplus_{i<j}\V^{(i)}\otimes\V^{(j)}$.
  We shall show that $\mathcal{Y}\subseteq\mathcal{U}$.
  Let $\mathcal{W}=\ker(\phi)$ be as in Remark~\ref{R:Spdq}
  where $\Lambda^2(\V)/\mathcal{W}$ is a 1-dimensional $EA$-module.
  If $i<j$, then $\V^{(i)}\otimes\V^{(j)}$ is a faithful $A^\infty$-module
  which is irreducible if $j\ne d/2+i$ by~\cite[\S5.4]{KL}, and is a sum of
  (two isomorphic) faithful irreducible submodules if $j=d/2+i$.
  Also, the uniserial proper $A^\infty$-submodule $\W$ in~Figure~\ref{F:Spdq} of
  $\Lambda^2(\V^{(i)})$ of dimension $2\ell^2-\ell-1$ is nontrivial.
  Since $A^\infty$ acts trivially on the simple factor module
  $A^2(V\otimes E)/\mathcal{U}$, we see that $\mathcal{U}$
  contains $\mathcal{Z}\coloneq\bigoplus_i\mathcal{W}^{(i)}\oplus\mathcal{Y}$,
  as claimed.

  Suppose now that $U$ is a maximal $F\kern-2pt A$-submodule of
  $A^2(V\otimes E)\restrict F$ such that $A^\infty$ acts trivially on
  $(A^2(V\otimes E)\restrict F)/U$. Choose 
  a maximal $EA$-submodule $\mathcal{U}$ of $A^2(V\otimes E)$ where $U$ contains
  $\mathcal{U}\restrict F$. Then $A^\infty$ acts trivially on
  $A^2(V\otimes E)/\mathcal{U}$. By the previous paragraph,
  $A^2(V\otimes E)/\mathcal{U}$ is a factor of
  $A^2(V\otimes E)/\mathcal{Z}=\bigoplus_i(\Lambda^2(\V)/\mathcal{W})^{(i)}$.
  Now $(A^2(V\otimes E)\restrict F)/U$ is an irreducible $F\kern-2pt A$-module
  and a factor $F\kern-2pt A$-module of
  $(A^2(V\otimes E)/\mathcal{U})\restrict F\cong
  (A^2(V\otimes E)\restrict F)/(\mathcal{U}\restrict F)$.
  The irreducible factor $F\kern-2pt A$-modules are isomorphic to
  a subfield $\F_{p^n}$ of $E=\F_{p^b}$ for some divisor $n$ of $b$
  by~\cite[VII.1.16(e)]{HB}. As $A$ varies, any divisor
  $n$ of $b$ can arise, see Lemma~\ref{L:SpTr}.\qed
\end{remark}

\begin{theorem}\label{T:Binftyeq1}
  Let $G$ be a finite nonabelian $3$-orbit $p$-group
  and let $V$, $A$, $W$, $B$ be as in Hypothesis~$\ref{H:pgp}$.
  If $A^\infty\ne1$ and  $B^\infty=1$,
  then $p$ is odd and $G,V,A,W,B$ are as in line $7$ of
  Table~$\ref{T:GVAWB}$ with $|\Phi(G)|=p^n$ as described in
  Lemma~$\ref{L:SpTr}$.
\end{theorem}

\begin{proof}
  If $p=2$, then $\Aut(G)$ is solvable by Theorem~\ref{T:Ainftyeq1} and so
  $A^\infty=1$, a contradiction.  Hence $p>2$. 
  If $n=1$, then $m$ must be even and $G$ is the extraspecial group of order
  $p^{1+m}$, and exponent~$p$ which appears on line~7 of Table~\ref{T:GVAWB}
  with $b=1$.
  Since $p>2$, we have $n\le\binom{m}{2}$ by Lemma~\ref{L:CentOdd}(a).
  Hence $m=2$ implies $n=1$. Suppose now that $m\ge3$ and $n\ge2$. Since
  $A^\infty\ne1$ and $B^\infty=1$, we have $A^\infty\le K$.
  If $H=N_{\GL(V)}(A^\infty)$, then $H/A^\infty\ge A/A^\infty\ge A/K\cong B$.
  We argue using Theorem~\ref{T:HeringNew} that $A^\infty\cong\Sp_{m/b}(p^b)$.
  Since $p^n-1$ divides $|B|$, we see that $B\ne1$. Hence $A$ is strictly
  larger than $A^\infty$, so case~(d) of Theorem~\ref{T:HeringNew} cannot
  hold, nor can case~(a) as $A^\infty\ne1$.
  In case~(c), we have $A^\infty=\SL_2(5)$ and $V=\F_3^{\,4}$ is an
  $A^\infty$-module. The maximal $A^\infty$-submodules of $\Lambda^2V$ have
  codimension 1 and $\binom{4}{2}-1=5$ by Remark~\ref{R:Spdq}.
  Therefore $n=1,5$ by
  Lemma~\ref{R:maximal}. But $n\ne1$, so $n=5$ and $242=p^n-1$ divides $|B|$.
  A direct calculation with {\sc Magma}~\cite{Magma} shows that $|H/A^\infty|=8$.
  This rules out case~(c). Case~(e) is similarly excluded. Hence we have
  $A^\infty\in\{\SL_{m/b}(p^b),\Sp_{m/b}(p^b)\}$. Suppose $A^\infty\cong\SL_{m/b}(p^b)$
  and $\V=\F_{p^b}^{\,m/b}$ is its natural module. As $\V$ is irreducible,
  so too is $\Lambda^2\V$ by~\eqref{E:SL3}. Let $V=\V\restrict\F_p=\F_p^{\,m}$.
  We claim that
  $\Lambda^2V$ is is a direct sum of faithful irreducible $\F_pA^\infty$-modules.
  The claim follows from Remark~\ref{R:C} as
  $\Lambda^2(\V\otimes\F_{p^b})=\bigoplus_i\Lambda^2(\V)^{(i)}
  \oplus\bigoplus_{i<j}\V^{(i)}\otimes\V^{(j)}$ and $\Lambda^2(\V)^{(i)}$
  is faithful and irreducible by~\eqref{E:SL3},
  and $\V^{(i)}\otimes\V^{(j)}$ is either a faithful and
  irreducible $\SL_{m/b}(p^b)$-modules or a direct sum of two such
  by~\cite[Theorem~5.4.5]{KL}.
  This implies that  $B^\infty\cong\SL_{m/b}(p^b)\ne1$, a contradiction.
  The only remaining possibility in Theorem~\ref{T:HeringNew}
  is $A^\infty\cong\Sp_{m/b}(p^b)$ where $m/b\ge4$ is even. In this case,
  $A\le\Gamma\Sp_{m/b}(p^b)$ since the normalizer
  of $\Sp_{m/b}(p^b)$ in $\GL_m(p)$ is $\Gamma\Sp_{m/b}(p^b)$. In summary,
  we have shown that
  $\Sp_{m/b}(p^b)=A^\infty\le A\lhdeq\Gamma\Sp_{m/b}(p^b)\le\GL_m(p)$.

  We now apply Theorem~\ref{T:HeringNew}(a) to $B\le\GL_n(p)$.
  First, $B^\infty=1$ and $B=A/K$ implies $\Sp_{m/b}(p^b)=A^\infty\le K$.
  Thus $B$ is a
  section of $\Gamma\Sp_{m/b}(p^b)/\Sp_{m/b}(p^b)\cong\GammaL_1(p^b)$, and
  so $B$ is metacyclic. Since $B^\infty=1$, the choices for $B$ are constrained
  by Theorem~\ref{T:HeringNew}(a).  The
  extraspecial group $D_8\circ Q_8=2_{-}^{1+4}$ is not metacyclic, and therefore
  $(n,p)\ne(4,3)$, as subgroups of metacyclic groups are metacyclic.
  Suppose that $n=2$ and $p\in\{5,7,11,23\}$.
  A calculation using {\sc Magma}~\cite{Magma} shows that the subgroups of
  $\GL_2(p)$ with $p\in\{5,7,11,23\}$ that are both metacyclic and
  transitive on nonzero vectors, all lie in $\GammaL_1(p^2)$.
  Therefore $B\le\GammaL_1(p^n)$ as in line~7 of Table~\ref{T:GVAWB}.

  By Lemma~\ref{L:CentOdd}(a) the 3-orbit
  group $G$ is isomorphic to $(V\rtimes\Lambda^2(V))/U$  where
  $V=\F_p^{\,m}$ is the natural $A$-module, and $U$ is a maximal submodule
  of $\Lambda^2(V)$ by Remark~\ref{R:maximal}.
  The simple quotient
  $A$-modules $\Lambda^2(V)/U$ of $\Lambda^2(V)$ are the
  subfields of $\F_{p^b}$ by Remark~\ref{R:sub}, and each subfield gives rise to
  a 3-orbit group $G$.
  Thus $G$ is as described on line~7 of Table~\ref{T:GVAWB}. Large
  subgroups of $A$ and $B$ are described in Lemma~\ref{L:SpTr}.
  Indeed, $A\le\Gamma\Sp_{m/b}(p^b)$ and $B\le\GammaL_1(p^b)$ as
  in line~7 of Table~\ref{T:GVAWB}.
\end{proof}

\section{Examples of \texorpdfstring{$k$}{}-orbit groups}\label{S:Ex}

In this section we give examples of $k$-orbit groups for small $k$. We
focus on 3-orbit groups. Extraspecial $p$-groups provide examples of both
3-orbit and 4-orbit groups.

  If $G$ is a finite extraspecial $p$-group, or an infinite Heisenberg group,
  then viewing the elements of $G$ as ordered pairs facilitates a
  geometric method to construct $\Aut(G)$. This method, was not used
  by Winter in~\cite{W}, but is used in Lemma~\ref{L:SpTr} below.
  
  \begin{remark}\label{R:Sp}
  We first describe $\GammaSp(\V)$ and $\CSp(\V)$.
  Let $f\colon\V\times\V\to F$ be a non-degenerate symplectic bilinear form
  on $\V=F^d$ where $d\ge2$ is even.
  Let $\GammaSp(\V)$ be the group of bijective semilinear symplectic
  similarities on $\V$.
  These satisfy
  \[
  (\lambda v)g =\lambda^{\sigma(g)}(vg), (v_1+v_2)g=v_1g+v_2g,
  \ \ \textup{and $f(v_1g,v_2g)=\delta(g)^{\sigma(g)}f(v_1,v_2)^{\sigma(g)}$,}
  \]
  for $g\in\GammaSp(\V)$, $v,v_1,v_2\in\V$,  $\lambda\in F$, where $\delta(g)\in F^\times$
  and $\sigma(g)\in\Aut(F)$ depend on $g$. The map
  $\sigma\colon\GammaSp(\V)\to\Aut(F)$ is an epimorphism.
  Comparing $f(v_1(gh),v_2(gh))$ to $f((v_1g)h,(v_2g)h)$ gives
  the the 1-cocycle condition
  $\delta(gh)=\delta(g)\delta(h)^{\sigma(g^{-1})}$.
  The \emph{conformal symplectic group} denoted by $\CSp(\V)$ is
  the kernel of $\sigma$ \emph{c.f.}~\cite[Def.~1.6.14]{BHRD}.

  We view the elements of $\GammaL_1(F)$ as products $\delta\sigma$
  with $(\delta,\sigma)\in F^\times\times\Aut(F)$ and multiplication rule
  $(\delta_1\sigma_1)(\delta_2\sigma_2)
  =\delta_1\delta_2^{\sigma_1^{-1}}\sigma_1\sigma_2$.
  It follows from the previous paragraph that
  $\GammaSp_d(F)=\Sp_d(F)\rtimes\GammaL_1(F)$, see~\cite[Table~2.1.C]{KL} and
  Remark~\ref{R:Spdq}.\qed
\end{remark}

 If $F:F_0$ is a Galois extension of the subfield $F_0$, then
 $\GammaL_1(F_0)$ is a factor group of $\GammaL_1(F)$, and hence
 $\GammaL_1(F_0)$ is a factor group of
 $\GammaSp_d(F)\cong\Sp(\V)\rtimes\GammaL_1(F)$.

\begin{lemma}\label{L:SpTr}
  Let $F:F_0$ be a finite Galois field extension where $\Char(F)=p\ge0$.
  Let $f\colon \V\times\V\to F$ be a
  non-degenerate alternating $F$-bilinear form on $\V=F^d$
  where $d$ is even. Let $\Tr$ be the trace map
  $F\to F_0\colon\lambda\mapsto\sum_{\sigma\in\Gal(F:F_0)}\lambda^{\sigma}$. 
  The set $G=\V\times F_0$ with the multiplication rule
  $(v_1,\zeta_1)(v_2,\zeta_2)=(v_1+v_2,\zeta_1+\zeta_2+\Tr(f(v_1,v_2)))$
  defines a group. 
  If $p\ne2$, then $G=G_{f,F_0}$ is a $3$-orbit group and
  \[
  \Sp_d(F)\rtimes(F_0^\times\rtimes\Aut(F))\le \Aut(G)^{G/G'}
  \qquad\textup{and}\qquad \GL_1(F_0)\le \Aut(G)\restrict G'.
  \]
  If $F=F_0=\F_p$, where $p$ is an odd prime, then $A=\CSp_d(p)$ and
  $B=\GL_1(p)$.
\end{lemma}

\begin{proof}
  Since the map
  $\V\times\V\to F_0\colon(v_1,v_2)\mapsto\Tr(f(v_1,v_2))$
  is biadditive, the multiplication on $G$ is associative.
  Hence $G$ is a group where $(0,0)$ is the identity element and
  $(v,\zeta)^{-1}=(-v,-\zeta)$ as $f(v,v)=0$. If $p>0$, then the
  exponent of $G$ is $p$ since
  $(v,\zeta)^k=(kv,k\zeta)$ for $k\in\mathbb{Z}$.
  The commutator $[(v_1,\zeta_1),(v_2,\zeta_2)]$
  equals $(0,2\Tr(f(v_1,v_2))$. Thus $G$ is abelian if $p=2$.
  Suppose now that $p\ne2$. As $f$ and $\Tr$ are surjective functions, 
  it follows that $G'=\{0\}\times F_0$.

  Let $\A$ be the subgroup of $\GammaSp(\V)$ (see Remark~\ref{R:Sp}) comprising
  all $g$ satisfying
  $f(v_1g,v_2g)=\delta(g)f(v_1,v_2)^{\sigma(g)}$ with
  $\delta(g)\in F_0^\times$. Then the structure of $\A$ is
  $\Sp_d(F)\rtimes(F_0^\times\rtimes\Aut(F))$. Using
  $\Tr(\delta\lambda^\sigma)=\delta\Tr(\lambda)$ for
  $\delta\in F_0$, $\lambda\in F$, $\sigma\in\Aut(F)$, we show below
  that $(v,\zeta)^g=(vg,\delta(g)\zeta)$ defines an action of
  $g\in\A$ on~$G$:
  \begin{align*}
    (v_1,\zeta_1)^g(v_2,\zeta_2)^g
    &=(v_1g,\delta(g)\zeta_1)(v_2g,\delta(g)\zeta_2)\\
    &=(v_1g+v_2g,
    \delta(g)\zeta_1+\delta(g)\zeta_2+\Tr(\delta(g)f(v_1,v_2)^{\sigma(g)}))\\
    &=((v_1+v_2)g,\delta(g)(\zeta_1+\zeta_2+\Tr(f(v_1,v_2))))\\
    &=(v_1+v_2,\zeta_1+\zeta_2+\Tr(f(v_1,v_2))^g=((v_1,\zeta_1)(v_2,\zeta_2))^g.
  \end{align*}
  Thus $g$ is a bijective endomorphism of $G$, i.e. an
  automorphism of $G$. Moreover, $\A$ acts on $G$ since
  $((v,\zeta)^g)^h=(v,\zeta)^{gh}$. Therefore $\Aut(G)$ has 3 orbits on~$G$,
  namely  $\{(0,0)\}$, $\{0\}\times F_0^\times$, $(\V\setminus\{0\})\times F_0$,
  that is $1, G'\setminus\{1\}, G\setminus G'$. We have therefore shown that
  $\A\le \Aut(G)^{G/G'}$ and $\GL_1(F_0)\le \Aut(G)\restrict G'$, as claimed.

  Finally, suppose that $F=F_0=\F_p$, where $p$ is an odd prime. In this
  case $G$ is an extraspecial $p$-group of order $p^{1+d}$ and exponent $p$. It
  follows from~\cite{W} that ${\rm Out}(p^{1+d})=\CSp_d(p)$ and hence
  $A=\CSp_d(p)$ and $B=\GL_1(p)$ as claimed.
\end{proof}

\begin{remark}\label{R:U}
  The subgroup $U$ in Remark~\ref{R:maximal} is the kernel of
  the map $\Lambda^2V\to F_0$ defined by
  $v_1\wedge v_2\mapsto\Tr(f(v_1,v_2))$
  where $f$ and $\Tr$ are as in Lemma~\ref{L:SpTr}.\qed
\end{remark}

  We prove that finite groups $G_{f,F_0}$ in Lemma~$\ref{L:SpTr}$
  of the same order are isomorphic. This will show that on line~7 of
  Table~\ref{T:GVAWB} there is only one isomorphism type of 3-orbit
  $p$-group $G$ with $p>2$, $|G/\Z(G)|=p^m$, $|\Z(G)|=p^n$ for each
  $n\mid m$ with $m/n$ even.
  
\begin{lemma}\label{L:GfGf'}
  Let $q$ be an odd prime power, and let $d,e$ be positive integers where $d$
  is even. Let $F=\F_{q^e}$, $\V=F^d$, $F'=\F_{q^{de/2}}$, $\V'=(F')^2$, and
  let $f'\colon\V'\times\V'\to F'$ be a non-degenerate alternating
  $F'$-bilinear form. Let $\phi\colon F'\to F^{d/2}$ be an $F$-isomorphism, and
  let $\Phi\colon\V'\to\V$ be an $F$-isomorphism obtained by applying $\phi$
  coordinatewise. Then the map $f\colon\V\times\V\to F$ defined by
  $f(v_1,v_2)=\Tr_{F'/F}(f'(\Phi^{-1}(v_1),\Phi^{-1}(v_2)))$ for $v_1,v_2\in\V$
  is a non-degenerate alternating $F$-bilinear form. If $F_0=\F_q$, then the
  groups $G_{f',F_0}$ and $G_{f,F_0}$ defined in Lemma~$\ref{L:SpTr}$
  are isomorphic of order $q^{de+1}$. 
\end{lemma}

\begin{proof}
  Clearly $f(v_2,v_1)=-f(v_1,v_2)$ and $f\colon\V\times\V\to F$ is an
  alternating $F$-bilinear form. Suppose $v_1\in\V$ is
  nonzero and $f(v_1,v_2)=0$ for all $v_2\in\V$, then $\Phi^{-1}(v_1)\ne0$ and
  $\Tr_{F'/F}(f'(\Phi^{-1}(v_1),\Phi^{-1}(v_2)))=0$
  implies $f'(\Phi^{-1}(v_1),\Phi^{-1}(v_2))=0$. As $f'$ is non-degenerate, we see
  that $\Phi^{-1}(v_1)=0$ and hence $v_1=0$, a contradiction. Therefore $f$ is
  non-degenerate. We define a group $H_{f'}=\V'\times F'$ with
  multiplication rule
  $(v'_1,z'_1)(v'_2,z'_2)=(v'_1+v'_2,z'_1+z'_2+f'(v'_1,v'_2))$. Similarly,
  let $H_f=\V\times F$ have multiplication rule
  $(v_1,z_1)(v_2,z_2)=(v_1+v_2,z_1+z_2+f(v_1,v_2))$.
  We show that
  $\psi\colon H_{f'}\to H_f$ defined by $(v',z')^\psi=(\Phi(v'),\Tr_{F'/F}(z'))$
  is an epimorphism. Certainly $\Tr_{F'/F}\colon F'\to F$ is an epimorphism, and
  $((v'_1,z'_1)(v'_2,z'_2))^\psi=(v'_1,z'_1)^\psi(v'_2,z'_2)^\psi$ holds because
  $\Tr_{F'/F}(f'(v'_1,v'_2))=f(\Phi(v'_1),\Phi(v'_2))$ for all $v'_1,v'_2\in\V'$.

  The epimorphism
  \[H_{f'}\to H_f/\ker(\Tr_{F/F_0})\quad\textup{defined by}\quad (v',z')\mapsto(\Phi(v'),\Tr_{F/F_0}(\Tr_{F'/F}(z')))\]
  maps $(v',z')$ to $(\Phi(v'),\Tr_{F'/F_0}(z'))$
  and has kernel $\ker(\Tr_{F'/F_0})$.
    The epimorphism $H_{f'}/\ker(\Tr_{F'/F_0})\to H_f/\ker(\Tr_{F/F_0})$ is
    the same as $\Psi\colon G_{f',F_0}\to G_{f,F_0}$. However,
    $|\V'|=|\V|=q^{de}$, $|F_0|=q$ and
    $|G_{f',F_0}|=q^{de+1}=|G_{f,F_0}|$, so
    $\Psi$ is an isomorphism.
\end{proof}

By Lemma~\ref{L:GfGf'} the groups $G_{f,\F_q}$ with $de$ constant are
isomorphic 3-orbit groups. By choosing a different factorization of
$de$ we see different subgroups of the automorphism group of
$G_{f,F_0}$ by Lemma~\ref{L:SpTr}.  Changing to the notation of line~7
of Table~$\ref{T:GVAWB}$, note that if $n\mid b\mid m$ where $m/b$ is even, then
$\Sp_{m/b}(p^b)\le\Sp_{m/n}(p^n)$ and both groups act transitively on a
set of $p^m-1$ nonzero vectors. Hence there is one isomorphism type
of 3-orbit group on line~7, namely with $b=n$ and
$A^\infty=\Sp_{m/n}(p^n)$.

Let $n_2$ denote the largest 2-power divisor of $n$, and let $n_{2'}=n/n_2$
be its odd~part.

\begin{theorem}\label{T:irred}
  An irredundant list of pairwise non-isomorphic $3$-orbit groups can be
  obtained from Table~$\ref{T:GVAWB}$ if on line~$3$, $G=A(n,\phi^i)$
  where $1\le i\le (n_{2'}-1)/2$ and $\phi=\theta^{n_2}$ has odd order
  $n_{2'}$ and $\theta(x)=x^2$ for $x\in\F_q$; and
  $b=n$ holds on line~$7$.
\end{theorem}

\begin{proof}
  The abelian groups $(\Cyc_{p^2})^n$ in line 1 of Table~\ref{T:GVAWB}
  are pairwise non-isomorphic (for different $n$), and not isomorphic
  to the (nonabelian) groups on other lines. Similarly, the groups $G$
  on line~2 have $|G/G'|=r$ and $|G'|=p^n$ and so are pairwise
  nonisomorphic for all prime divisors $r$ of $n-1$, and not isomorphic
  to the groups on other lines. The only isomorphisms between the
  line 3 groups are $A(n,\theta)\cong A(n,\theta^{-1})$
  by~\cite[Theorem~2]{Hig}. Hence $\phi=\theta^{n_2}$
  has odd order $n/n_2=n_{2'}$, and the isomorphism classes
  are $A(n,\phi^i)$ for $1\le i\le (n_{2'}-1)/2$.
  Comparing $|G|$, $|\Z(G)|$ and $|A|$, there are no other isomorphisms
  between the 2-groups on lines 3--5 since $B(3)\not\cong P$.
  By Lemma~\ref{L:GfGf'} and the above remarks, the group on
  line~7 has $b=n$ and
  $A^\infty=\Sp_{m/n}(p^n)$ where $m/n$ is even. Since
  $n\le m/2<m$, the group $G$ on line~7 has $p^n=|\Z(G)|<|G/\Z(G)|=p^m$ and so
  is not isomorphic to the group on line~6. This is completes the proof.
\end{proof}

\begin{remark}\label{R:irred}
    A major contribution of~\cite{LZ2} is to classify the odd order
    3-orbit $p$-groups~\cite[lines 6-8, Table~1]{LZ2}. The authors sought
    a complete list, but not an irredundant one.
    The orbits of $\GL(V)$ on subspaces of $\Lambda^2(V)$ of a fixed codimension
    are poorly understood, and this is at the heart of the difficulty of the
    $p$-group isomorphism problem. Indeed, the group isomorphism problem
    in~\cite[Lemma 2.12]{LZ2} can be solved if the (difficult) orbit problem
    on $\Lambda^2(V)$ can be solved. The approach here is to consider the
    group $\mathcal{G}=\V\times\Lambda^2(\V)$ with multiplication
    rule $(v_1,w_1)(v_2,w_2)=(v_1+v_2,w_1+w_2+v_1\wedge v_2)$ on which
    $\Gamma\Sp(\V)$ acts, and the groups on line~7 of Table~1 are factor groups
    of $\mathcal{G}$.   The approach in~\cite{LZ2} is to use a finitely
    presented group $\mathcal{H}_{n,p}$~\cite[Definition~2.9]{LZ2},
    or a certain central product $q_{+}^{1+m}$~\cite[Definition~2.17]{LZ2},
    with $m$ even,  and factor out a central subgroup $U$
    satisfying~\cite[Example~5.3(i)-(ii)]{LZ2}. Such a $U$ contains
    no ``subfield hyperplanes''~\cite[Definition~5.1]{LZ2}. The approach here
    involves considering quotient groups $\mathcal{G}/U$ where $U$ is a kernel
    of a trace map, and $U$ \emph{is} a subfield hyperplane. The fact that
    trace maps satisfy $\Tr_{K/F}\circ\Tr_{L/K}=\Tr_{L/F}$ helps to prove Lemma~\ref{L:GfGf'} and hence to give the irredundant list in Theorem~\ref{T:irred}.\qed
\end{remark}

\begin{lemma}
  The group $B(n)$ in Definition~{\rm \ref{D:ABP}(b)} is isomorphic to the
  Suzuki $2$-group $B(n,1,\xi)$ defined by~{\rm\cite[{\it Column~III}\,]{Hig}}
  where $\xi\ne\tau+\tau^{-1}$ for all $\tau\in \F_{2^n}^\times$.
\end{lemma}

\begin{proof}
  Let $q=2^n$.
  The polynomial $t^2+\xi t+1$ is irreducible in $\F_q[t]$ since
  $\tau+\tau^{-1}=\xi$ has no solutions for $\tau\in\F_q^\times$. Let $\eps$ be a
  root of $t^2+\xi t+1$. Then $\eps+\eps^q=\xi$ and $\eps^{q+1}=1$.
  Hence $\F_q[\eps]=\F_{q^2}$ and the norm map $\F_{q^2}^\times\to\F_q^\times$ sends
  $\alpha+\beta\eps\in \F_{q^2}^\times$ to $(\alpha+\beta\eps)(\alpha+\beta\eps^q)=\alpha^2+\xi \alpha\beta+\beta^2$. The Suzuki 2-group $B(n,1,\xi)$ can,
  by Higman~\cite[Column~V]{Hig}, be identified
  with the set $\F_q^{\,3}$ with the following multiplication rule
  \[
    (\alpha_1,\beta_1,\zeta_1)(\alpha_2,\beta_2,\zeta_2)=(\alpha_1+\alpha_2,\beta_1+\beta_2,\zeta_1+\zeta_2+\alpha_1\alpha_2+\xi \alpha_1\beta_2+\beta_1\beta_2).
  \]
  The third coordinate is related to the `bilinearized' form of the norm map
  \[
  (\alpha_1+\beta_1\eps)(\alpha_2+\beta_2\eps)^q=(\alpha_1+\beta_1\eps)(\alpha_2+\beta_2\eps^{-1})=\alpha_1\alpha_2+\xi \alpha_1\beta_2+\beta_1\beta_2.
  \]
  Therefore, $(\alpha_1+\beta_1\eps)(\alpha_2+\beta_2\eps)^q\in\F_q$ and hence
  \[
    (\alpha_1+\beta_1\eps)(\alpha_2+\beta_2\eps)^q\eps+((\alpha_1+\beta_1\eps)(\alpha_2+\beta_2\eps)^q\eps)^q=(\alpha_1\alpha_2+\xi \alpha_1\beta_2+\beta_1\beta_2)(\eps+\eps^q).
  \]
  Since $\eps+\eps^q=\xi$, the map
  $B(n,1,\xi)\to B(n)$ defined by
  $(\alpha,\beta,\zeta)\mapsto(\alpha+\beta\eps,\zeta\xi)$
  is an isomorphism. Consequently, the isomorphism type of $B(n,1,\xi)$ is
  independent of the choice of $\xi$ for which $t^2+\xi t+1$ is irreducible.
\end{proof}

We will construct examples of 3- and 4-orbit groups using the exterior
algebra $\Lambda(\V)$ of a vector space $\V$. If $\dim(\V)=d$, then
$\Lambda(\V)=\bigoplus_{k=0}^{\dim(V)}\Lambda^k(\V)$ is a graded algebra with $\dim(\Lambda^k(\V))=\binom{d}{k}$ and hence $\dim(\Lambda(\V))=2^d$. 
The following preliminary lemma exploits the action of $\GL(\V)$ 
on $\Lambda^k(\V)$, see~\cite[XIX]{Lang}.

\begin{lemma}\label{L:AltKer}
  Let $\V=F^d$ be an $d$-dimensional vector space over a field $F$.
  Suppose that $1< k\le d$ and $n=\binom{d}{k}$. The action of $\GL(\V)$
  on $\Lambda^k(\V)$ induces a homomorphism
  $\phi_{d,k}\colon\GL_d(F)\to\GL_n(F)$ of matrix groups.
  The kernel of $\phi_{d,k}$ is $\GL_d(F)$ if $k>d$, $\SL_d(F)$ if $k=d$, and
  $\{\lambda I_d\mid \lambda\in F\ {\rm and}\ \lambda^k=1\}$ if $k<d$.
\end{lemma}

\begin{proof}
  If $k>d$, then $\Lambda^k(\V)=\{0\}$ so $\ker\phi_{d,k}=\GL_d(F)$.
  Let $\V=\langle e_1,\dots, e_d\rangle$.
  If $k=d$, then $\Lambda^d(\V)=\langle e_1\wedge\cdots\wedge e_d\rangle$ and
  $g\phi_{d,d}=\begin{pmatrix}\det(g)\end{pmatrix}$, so that
  $\ker\phi_{d,d}=\SL_d(F)$.

  Let $\langle v_1,\dots,v_k\rangle$ be a typical $k$-subspace of $\V$ where $k<d$.
  As $g\in\ker\phi_{d,k}$ fixes $v_1\wedge\cdots\wedge v_k$, it also fixes
  the  $k$-subspace $\langle v_1,\dots,v_k\rangle$
  by~\cite[Lemma~12.6]{Taylor}. As $k<d$ we may choose a vector $v_{k+1}$ in
  $\V\setminus\langle v_1,\dots,v_k\rangle$. Since $g$ fixes the $k$-subspaces
  $\langle v_1,\dots,v_k\rangle$ and $\langle v_2,\dots,v_{k+1}\rangle$, it fixes
  their intersection, \emph{viz.} $\langle v_2,\dots,v_{k}\rangle$. Thus $g$
  fixes all $(k-1)$-subspaces. By induction, $g$ fixes all 1-subspaces of $\V$
  and hence $g$ is a scalar matrix. However, $\lambda I_d\in\ker\phi_{d,k}$
  precisely when $\lambda^k=1$. This completes the proof.
\end{proof}

\begin{remark}
  If $F=\F_q$, then $\{\lambda\in\F_q^\times\mid \lambda^k=1\}$
  is cyclic of order~$\gcd(k,q-1)$.\qed  
\end{remark}

\begin{lemma}\label{L:GammaL3}
  Let $\Lambda(\V)$ be the exterior algebra of the $F$-vector space $\V=F^3$
  where $\textup{char}(F)\ne2$. Then the set
  $G=G_F=\V\times\Lambda^2\V\times\Lambda^3\V$ with the multiplication~rule
  \[
  (v_1,w_1,x_1)(v_2,w_2,x_2)=
  (v_1+v_2,w_1+w_1+v_1\wedge v_2,x_1+x_2+v_1\wedge w_2+w_1\wedge v_2)
  \]
  defines a $4$-orbit group. Also $\Aut(G)$ induces on
  $G/\gamma_2(G)$, $\gamma_2(G)/\gamma_3(G)$ and $\gamma_3(G)$ subgroups
  $A$, $B$, $C$ respectively where
  $\GammaL(\V)\le A$,
  $\{g\wedge g\mid g\in\GammaL(\V)\}\le B$
  and $\{g\wedge g\wedge g\mid g\in\GammaL(\V)\}\le C$.
  In particular, $\gamma_3(G)=\Z(G)$, $\gamma_2(G)=C_G(\gamma_2(G))$ and $G/\gamma_3(G)$
  is a $3$-orbit group. If $|F|=q$ is odd,
  then $|G|=q^7$ and $|G/\gamma_3(G)|=q^6$.
\end{lemma}

\begin{proof}
  The exterior algebra $\Lambda(\V)$ equals
  $\bigoplus_{i=0}^3\Lambda^i(\V)$ where
  $\dim(\Lambda^i(\V))=\binom{3}{i}$.
  A basis $(e_1,e_2,e_3)$ for $\Lambda^1(\V)=\V$ gives bases
  $(e_2\wedge e_3,e_3\wedge e_1,e_1\wedge e_2)$
  for $\Lambda^2(\V)$ and $(e_1\wedge e_2\wedge e_3)$ for $\Lambda^3(\V)$.
  Relative to these bases a $3\times 3$ matrix $g\in\GL(\V)$ induces
  the $3\times 3$ matrix
  $\det(g)g^{-{\rm T}}=g\wedge g\in\GL(\Lambda^2\V)$ and the $1\times 1$ matrix
  $(\det(g))=g\wedge g\wedge g\in\GL(\Lambda^3\V)$. Hence the action of
  $\GL(\V)$ on $\Lambda^2(\V)$ is different from the `natural' and `dual' actions.

  The group of units $\Lambda(\V)^\times$ has a normal subgroup
  $M=\{1\}\times\V\times\Lambda^2\V\times\Lambda^3\V$, and
  \begin{align*}
  (&1+v_1+w_1+x_1)\wedge(1+v_2+w_2+x_2)\\
  =\;&1+(v_1+v_2)+(w_1+v_1\wedge v_2+w_2)+(x_1+w_1\wedge v_2+v_1\wedge w_2+x_2).
  \end{align*}
  Therefore the stated multiplication rule of triples
  in $G=\V\times\Lambda^2\V\times\Lambda^3\V$ defines an isomorphism
  $G\to M\colon (v,w,x)\mapsto 1+v+w+x$. In particular, $G$ is a group.

  The identity element of $G$ is $(0,0,0)$
  and $(v,w,x)^{-1}=(-v,-w,-x)$ since
  $\wedge$ is antisymmetric and $w\wedge v+v\wedge w=0$. Since
  $(v,w,x)^k=(kv,kw,kx)$ for $k\in\ZZ$ it follows that $G$ is 
  torsion free if $\Char(F)=0$, and has (odd) exponent $p=\Char(F)$ otherwise.
  In both cases $\gamma_2(G)=\{0\}\times\Lambda^2\V\times\Lambda^3\V$ holds
  because
  \begin{equation}\label{E:comm2}
    [(v_1,w_1,x_1),(v_2,w_2,x_2)]
    =(0,2v_1\wedge v_2,2(v_1\wedge w_2+w_1\wedge v_2)).
  \end{equation}
  Setting $[(v_1,w_1,x_1),(v_2,w_2,x_2)]=(0,w',x')$ in~\eqref{E:comm2}
  where $w'=2v_1\wedge v_2$ gives that
  \[
    [[(v_1,w_1,x_1),(v_2,w_2,x_2)],(v_3,w_3,x_3)]=(0,0,2w'\wedge v_3)
    =(0,0,4v_1\wedge v_2\wedge v_3).
  \]
  Hence $\gamma_3(G)=\{0\}\times\{0\}\times\Lambda^3\V$ as $\Char(F)\ne2$.
  Observe that if $v\in\V$ satisfies $v\wedge w=0$ for all $w\in\Lambda^2\V$,
  then $v=0$. Hence~\eqref{E:comm2} implies that $\gamma_3(G)=\Z(G)$,
  and $C_G(\gamma_2(G))=\gamma_2(G)$.  

  Now $g\in\GL(\V)$ acts on $G$ is via
  $(v,w,x)^g=(vg,w(g\wedge g),x(g\wedge g\wedge g))$ as described above.
  Hence $G$ is transitive on the nonzero vectors of $G/\gamma_2(G)=V$,
  $\gamma_2(G)/\gamma_3(G)=\Lambda^2V$, $\gamma_3(G)/\gamma_4(G)=\Lambda^3V$,
  so $G$ is a 4-orbit group.
  Further, $\sigma\in\Aut(F)$ acts to $G$ via
  $(v,w,x)^\sigma=(v^\sigma,w^\sigma,x^\sigma)$ by applying $\sigma$ to the
  coordinates of $v,w,x$ relative the stated bases. This shows that
  $\GammaL(\V)$ is a subgroup of $\Aut(G)$ which induces
  $\GammaL(\V)\le A$, $\{g\wedge g\mid g\in\GammaL(\V)\}\le B$
  and $\GammaL_1(F)=\{g\wedge g\wedge g\mid g\in\GammaL(\V)\}\le C$ as claimed.
\end{proof}

\begin{remark}\label{R:det}
  The group $G_F$ in Lemma~\ref{L:GammaL3} is abelian if $\Char(F)=2$.
  If $g\in\GL_3(F)$, then $g\wedge g=\det(g)g^{-{\rm T}}$ so that
  $\det(g\wedge g)=\det(g)^3\det(g^{-{\rm T}})=\det(g)^2\in(F^\times)^2$.
  The homomorphism $\GL(\V)\to\GL(\Lambda^2\V)$ has kernel $\langle-1\rangle$,
  and $\GL(\V)\to\GL(\Lambda^3\V)$ has kernel $\SL(\V)$
  by Lemma~\ref{L:AltKer}.\qed
\end{remark}


\begin{lemma}\label{L:GL3}
  Let $F:F_0$ be a finite separable field extension where $\Char(F)\ne2$.
  Then the trace map $\Tr\colon F\to F_0$ is surjective, and 
  the set $G=G_{F,F_0}=F^3\times F^3\times F_0$ endowed with
  the multiplication rule
  \[
  (v_1,w_1,x_1)(v_2,w_2,x_2)=
  (v_1+v_2,w_1+w_1+v_1\wedge v_2,x_1+x_2+\Tr(v_1\wedge w_2+w_1\wedge v_2))
  \]
  defines a group.
  Let $H=\{g\in\GL_3(F)\mid \det(g)\in F_0^\times\}$ and let $H^{+}$ be the
  subgroup $H^{+}=\{g\in\GL_3(F)\mid \det(g)\in (F^\times)^2\}$.
  If $\Aut(G)$ induces on
  $G/\gamma_2(G)$, $\gamma_2(G)/\gamma_3(G)$ and $\gamma_3(G)$ subgroups
  $A$, $B$, $C$ respectively, then
  $H\rtimes\Aut(F)\le A$, $H^{+}\rtimes\Aut(F)\le B$ and $F_0^\times\le C$.
  Moreover, $G$ is a $4$-orbit group and $G/\gamma_3(G)$ is
  a $3$-orbit group.
\end{lemma}

\begin{proof}
  If $\sigma_1,\dots,\sigma_{|F:F_0|}$ are the $F_0$-linear embeddings
  $F\to\overline{F}$ into the algebraic closure $\overline{F}$ of $F$, then
  $\Tr(x)=\sum_{i=1}^{|F:F_0|}\sigma_i(x)$. Since $\sigma_1,\dots,\sigma_n$ are
  linearly independent over $F_0$, the $F_0$-linear map $\Tr$ is
  nonzero, and hence is surjective.  Let $M$ be the kernel of the trace
  map $\Tr\colon F\to F_0$. The First Isomorphism Theorem
  gives $F^{+}/M\cong F_0^{+}$.
  We will show that $G_{F,F_0}$ is a factor group of the group $G_F$
  in Lemma~\ref{L:GammaL3}. Indeed, the map
  $\phi\colon G_F\to G_{F,F_0}\colon (v,w,x)\mapsto (v,w,x+M)$
  preserves multiplication and has kernel
  $\{(0,0,x)\mid x\in M\}\le\Z(G_F)$ where $x+M$
  is viewed as an element of~$F_0$ via the isomorphism
  $F^{+}/M\cong F_0^{+}$. Therefore $G_F/M\cong G_{F,F_0}$.

  The epimorphism $\phi$ maps $\gamma_i(G_F)$ to $\gamma_i(G_{F,F_0})$
  for $1\le i\le 3$. Further, if $g\in H$, then $g\wedge g\in H^{+}$
  by Remark~\ref{R:det}. The remaining claims follow from Lemma~\ref{L:GammaL3}.
\end{proof}

Subgroups $G_1,G_2\le\Sym(\Omega)$ with the same orbits
on $\Omega$ are called \emph{orbit-equivalent}.

\begin{lemma}\label{L:prb}
  Let $F$ be a division ring, and $C\le F^\times$ a finite subgroup.
  Suppose that $A\le\Aut(F)$ fixes $C$ setwise and is orbit equivalent to
  $\Aut(C)\le\Sym(C)$.
  Let $\V=F^d$ be a $d$-space over $F$.
  Then the set $G=C\times\V$ endowed with the multiplication rule
  $(\lambda,v)(\mu,w)=(\lambda\mu,\mu v+w)$
  defines a group. Further,
  $\Aut(G)$ has one more orbit on~$G$ than $\Aut(C)$ has on $C$,
  i.e. $\omega(G)=\omega(C)+1$.
\end{lemma}

\begin{proof}
  We now show that the set $G=C\times\V$ is a group.
  Associativity holds as
  \[
  ((\lambda,v)(\mu,w))(\nu,x)=(\lambda\mu\nu,\mu\nu v+\nu w+x)
  =(\lambda,v)((\mu,w)(\nu,x))
  \]
  holds for all $(\lambda,v),(\mu,w),(\nu,x)\in C\times \V$.
  The identity element of $G$ is $(1,0)$. Also $(\lambda,v)$ has inverse
  $(\lambda^{-1},-\lambda^{-1}v)$
  and $(\lambda,v)^n=(\lambda^n,(\lambda^{n-1}+\cdots+\lambda+1)v)$ for $n\ge0$.

  We now show that $M\coloneq\{1\}\times\V$ is characteristic in $G$.
  If $\Char(F)=0$, then this follows since elements of $G\setminus M$
  have finite order (as $\lambda^n=1$ implies $(\lambda,v)^n=(1,0)$
  if $\lambda\ne1$), while nontrivial elements of $M$ have infinite order.
  If $\Char(F)>0$, then $C$ is contained in the multiplicative
  group of a finite field by the proof of~\cite[Theorem~6]{Herstein}.
  Hence $M$ is a normal
  Sylow $p$-subgroup of $G$ and thus characteristic in $G$.
  We next show that
  $\omega(G)=\omega(C)+\omega(M)-1$.
  Clearly $G/M\cong C$. First, $\Aut(M)$ has two orbits on $M$.
  Note that an invertible $F$-linear map $g\in\Aut_F(\V)$ acts on $G$ via
  $(\mu,v)^g=(\mu,vg)$. Hence $\Aut_F(\V)$ has one nontrivial orbit on $M$,
  so $\omega(M)=2$.   Also $\alpha\in A\le\Aut(F)$ acts coordinatewise
  on $M\cong F^d$, and hence $A$
  acts on $G$ via $(\lambda,v)^\alpha=(\lambda^\alpha,v^\alpha)$. Since $A$
  is orbit-equivalent to $\Aut(C)\le\Sym(C)$, both $A$ and $\Aut(C)$ have
  $\omega(C)$ orbits on $C$.
  These two types of automorphisms of $G$ generate a subgroup of $\Aut(G)$
  with $\omega(C)+1$ orbits. Hence $\omega(G)=\omega(C)+1$
  by Lemma~\ref{L:C}.
\end{proof}

\begin{example}
  Let $F=\mathbb{H}$ be the real quaternions. Then $\mathbb{H}^\times$
  contains the quaternion subgroup $C=\{\pm1,\pm i,\pm j,\pm k\}\cong Q_8$.
  Set $r=\frac{i+j}{\sqrt2}$, $s=\frac{1+i+j+k}{2}$ and $t=\frac{1+i}{\sqrt2}$.
  The binary octahedral group
  $\textup{BO}=\langle r,s,t\mid r^2=s^3=t^4=rst\rangle$ satisfies
  $|\textup{BO}|=48$, $\Z(\textup{BO})=\langle rst\rangle=\langle -1\rangle$,
  and $C\lhdeq\textup{BO}$.
  The subgroup $A$ of $\Aut(F)$ comprising the inner automorphisms
  $F\to F\colon\lambda\mapsto\alpha^{-1}\lambda\alpha$, $\alpha\in\textup{BO}$,
  fixes $\textup{BO}'= Q_8=C$ setwise. Furthermore
  $A\cong\textup{BO}/\langle -1\rangle\cong S_4\cong\Aut(Q_8)$,
  so $A$ is orbit-equivalent to $\Aut(C)$. 
  Lemma~\ref{L:prb} shows that $G=C\times\mathbb{H}^d$ satisfies
  $\omega(G)=\omega(Q_8)+1=4$ for $d\ge1$.\qed
\end{example}

\begin{example}\label{E:Gpr}
  Let $p,r$ be distinct primes. Set $e=r^\ell$ where $\ell\ge1$.
  Then $p\nmid e$ and the cyclotomic
  polynomial $\Phi_e(t)$ is irreducible over the field $\F_p$
  precisely when $p$ has order $\deg(\Phi_e(t))=\phi(e)$ modulo $e$
  by~\cite[Ex.\;3.42,\;p.\,124]{LN}; in this case $\F_q$ where $q=p^{\phi(e)}$
  is the splitting field of $\Phi_e(t)$ over $\F_p$. Hence $\Cyc_e\le\F_q^\times$
  and $\Gal(\F_q/\F_p)=\Cyc_{\phi(e)}\cong\Aut(\Cyc_e)$. Set $C=\Cyc_e$, $F=\F_q$
  and $A=\Gal(\F_q/\F_p)$ in Lemma~\ref{L:prb} noting that the orbits
  of $\Aut(C)$ and $A$ on $C$ are the elements of $C$ of the same order.
  Thus the set $G=\Cyc_e\times\V$ is a group $\Cyc_e\ltimes\F_q^{\;d}$
  with $\omega(G)=\omega(\Cyc_{r^\ell})+1=\ell+2$. Setting $\ell=1$ gives
  $e=r$, $q=p^{r-1}$, $G=\Cyc_r\ltimes(\Cyc_{p})^{d(r-1)}$, and
  $\omega(G)=3$ as in line~2 of Table~\ref{T:GVAWB}.\qed
\end{example}

\begin{lemma}\label{L:prb2}
  If $p\ne r$ are prime and $G$ is an $(\ell+2)$-orbit group with
  $|G|=p^mr^\ell$, $G'\cong\Cyc_p^{\,m}$ and $G/G'\cong\Cyc_{r^\ell}$, then $G$ is
  isomorphic to the group in Example~$\ref{E:Gpr}$.
\end{lemma}

\begin{proof}
  By assumption, $\Aut(G)$ has precisely $\ell+2$ orbits on $G$.
  As $|\ord{G}|=\ell+2$, these orbits are the sets
  $O_1, O_p, O_r,\dots, O_{r^\ell}$ where $O_n=\{g\in G\mid |g|=n\}$.
  Let $R$ be a Sylow $r$-subgroup, and
  let $P=G'$ be the normal $p$-subgroup. Now $\Z(G)\cap P$ is trivial,
  otherwise $G$ has at least $\ell+3$ orbits. Hence $R$ acts fixed-point-freely
  on $P$. By Maschke's theorem $P=P_1\oplus\cdots\oplus P_d$ where each
  $P_i$ is an irreducible $R$-module. The $P_i$ must be
  pairwise isomorphic $R$-modules, otherwise $\Aut(G)$ has at
  least 3 orbits on $P$, a contradiction. Let $|P_1|=\cdots=|P_d|=p^b=q$.
  Hence each $\lambda\in R$ may be viewed as acting as
  a $d\times d$ scalar matrix over $\F_q$.
  Thus $G\cong\F_q^d\rtimes\Cyc_{r^\ell}$
  in Example~\ref{E:Gpr}.
\end{proof}

\begin{lemma}\label{L:abc}
  Let $q=2^n$ and let $\Tr\colon\F_{q^2}\to\F_q\colon\mu\mapsto\mu+\mu^q$ denote
  the trace map.
  \begin{enumerate}[{\rm (a)}]
  \item If $\theta\in\Gal(\F_q/\F_2)$ has $|\theta|>1$ odd, then $n\ne2^\ell$ and
    the group $A(n,\theta)$ in Definition~$\ref{D:ABP}$\textup{(a)} is a
    $3$-orbit $2$-group of order $q^2$ with $A\cong B\cong\GammaL_1(q)$.
  \item If $n\ge1$ and $\eps\in\F_{q^2}^\times$ have order $q+1$,
    then $B(n)=B_\eps(n)$ in Definition~$\ref{D:ABP}$\textup{(b)}
    is a  $3$-orbit $2$-group of order $q^3$ with $A\cong\GammaL_1(q^2)$
    and $B\cong\GammaL_1(q)$.
  \item Let $q=8$ and $\F_{q^2}^\times=\langle\eps\rangle\cong\Cyc_{63}$.
    The group $P=P(\eps)$  in Definition~$\ref{D:ABP}$\textup{(c)}
    is a $3$-orbit $2$-group isomorphic to
    \textup{SmallGroup($2^9$,10\,481\,182)} in {\sc Magma}~{\rm\cite{Magma}}
    with $A\cong\Cyc_7\rtimes\Cyc_9<\GammaL_1(\F_{64})$ of order $63$
    and $B=\GammaL_1(\F_8)$ of order $21$.
  \end{enumerate}
\end{lemma}

\begin{proof}
  (a)~A simple calculation shows that $G\coloneq A(n,\theta)$ defined by
  Definition~\ref{D:ABP}(a) is a group with
  $\Z(G)=G'=\Phi(G)=\{0\}\times\F_q$,
  see also~\cite[Theorem~(iv),\,p.\,704]{Dornhoff}.
  The group $\GammaL_1(q)$ acts faithfully on $G$ via
  $(\mu,\zeta)^{(\alpha,\lambda)}
  =(\mu^\alpha\lambda,\zeta^\alpha\lambda\lambda^\theta)$.
  Hence $G$ is a $3$-orbit $2$-group of order $q^2$ with $m=n$
  and $\GammaL_1(q)\le A$ and $\GammaL_1(q)\le B$.
  No solvable group of $\GL_n(2)$  properly contains $\GammaL_1(2^n)$
  by Theorem~\ref{T:HeringNew}(a). Therefore $A=\GammaL_1(2^n)$.
  It follows from~\cite{Hig} that $\GammaL_1(q)\le B$. Similar reasoning
  shows that  $B=\GammaL_1(2^n)$, so line~3 of Table~\ref{T:GVAWB} is valid.

  (b)~The group $B_\eps(n)$ appears in~\cite[Theorem~(v)]{Dornhoff}. Since
  $B_\eps(n)\cong B_{\eps'}(n)$ when
  $\langle\eps\rangle=\langle\eps'\rangle$ has order $q+1$ we
  write $B(n)$ instead of $B_\eps(n)$.
  The multiplication rule
  in~\cite[Theorem~(v)]{Dornhoff} can be rewritten as
  $(\mu_1,\zeta_1)(\mu_2,\zeta_2)
  =(\mu_1+\mu_2,\zeta_1+\zeta_2+\Tr(\eps\mu_1\mu_2^q))$
  since $\eps^q=\eps^{-1}$.
  The group $\GammaL_1(q^2)$ acts faithfully on $B(n)$ via
  $(\mu,\zeta)^{(\alpha,\lambda)}
  =(\mu^\alpha\lambda,\zeta^\alpha\lambda\lambda^\theta)$. Arguing as in part~(a)
  we have $A\cong\GammaL_1(q^2)$ and $B\cong\GammaL_1(q)$ as in
  line~4 of Table~\ref{T:GVAWB}.

  (c)~The group $P(\eps)$ appears in~\cite[Theorem~(vi)]{Dornhoff}.
  If $\F_{q^2}^\times=\langle\eps\rangle=\langle\eps'\rangle$, then
  $P(\eps)\cong P(\eps')$ so we write $P$ rather than $P(\eps)$.
  The maps $\psi,\phi\colon P\to P$ defined by
  \[
  (\alpha,\zeta)^\psi=(\eps^3\alpha,\eps^9\zeta)
  \qquad\textup{and}\qquad
  (\alpha,\zeta)^\phi=(\eps\alpha^4,\zeta^4)
  \]
  can be shown to be homomorphisms of $P$ that satisfy
  $\psi^{21}=\phi^9=1$, $\psi^\phi=\psi^4$ and $\psi^7=\phi^3$. Hence
  $\psi,\phi\in\Aut(P)$ and
  $\langle \psi,\phi\rangle=\langle \psi^3,\phi\rangle=\Cyc_7\rtimes\Cyc_9$.
  A computation with {\sc Magma}~\cite{Magma} shows that $|\Aut(P)|=2^{18}\cdot63$. There
  are $2^{18}$ central automorphisms so
  $A\cong\Cyc_{21}.\Cyc_3\cong\Cyc_7\rtimes\Cyc_9<\GammaL_1(\F_{64})$ and
  $B=\Cyc_7\rtimes\Cyc_3=\GammaL_1(\F_8)$ as
  $\F_8^\times=\langle\eps^9\rangle$.
\end{proof}

\section{Examples of 4-orbit groups}\label{S:4orbit}

\begin{center}
\begin{minipage}[b]{0.9\linewidth}
  \noindent
  ``A good stock of examples, as large as possible, is indispensable for
  a thorough understanding of any concept, and when I want to learn something
  new, I make it my first job to build one.''\hfill\sc{Paul Halmos}
\end{minipage}
\end{center}
  
In this section, we consider the feasibility of classifying finite $k$-orbit groups
for $k\le 6$. The nonsolvable $k$-orbit groups have been classified for
$k=4,5,6$, see \S1 for details. To assess the feasibility of classifying
the solvable $k$-orbit groups for $k=4,5,6$, we employ Halmos' strategy,
and seek a large stock of examples, particularly when $k=4$.
Using {\sc Magma}~\cite{Magma}, we studied the 1265679 groups of order less
than $2^{10}$ excluding $2^9$. Only 86 of these are nonabelian
solvable 4-orbit groups!
It appears that these groups belong to a small number of
infinite (and finite) families. For brevity, we list (without proof
of correctness) many of these below and in Table~\ref{T:4orbit}.
Given the difficulty of computing automorphism
groups, classifying the solvable 4-orbit groups may just be feasible. The most
difficult case will be when $G$ does not have four `obvious' characteristic
subsets (e.g. determined by element orders or characteristic subgroups).
For $k=5,6$ a complete classification may involve too many possibilities.

The $\Aut(G)$-orbit lengths for a 3-orbit group $G$ follow from
Theorem~\ref{T:main}. They are $1,p^n-1,p^n(r^m-1)$ where $p=r$ except for
line~2 of Table~\ref{T:GVAWB}, and the respective orbit-element orders
are $1,p,p^2$ in lines $1,3,4,5$ and $1,p,r$ otherwise.
For $k$-orbit groups with $k\ge4$ the orbit lengths and
orders are less obvious. Clearly the sum of the orbit lengths is $|G|$
and some orders in $\ord{G}$ may be be duplicated, see Table~\ref{T:4orbit}.
If $G$ is a solvable 4-orbit group with precisely 4 characteristic subgroups,
arranged as
\begin{tikzpicture}[scale=0.2]
\draw [-,semithick] (1,0) -- (0,1) -- (1,2) -- (2,1) -- (1,0);
\draw [fill=black] (1,0) circle [radius=1.4mm];
\draw [fill=black] ((0,1) circle [radius=1.4mm];
\draw [fill=black] ((1,2) circle [radius=1.4mm];
\draw [fill=black] ((2,1) circle [radius=1.4mm];
\end{tikzpicture}
or \,
\begin{tikzpicture}[scale=0.12]
\draw [-,semithick] (0,0) -- (0,3);
\draw [fill=black] (0,0) circle [radius=2.4mm];
\draw [fill=black] ((0,1) circle [radius=2.4mm];
\draw [fill=black] ((0,2) circle [radius=2.4mm];
\draw [fill=black] ((0,3) circle [radius=2.4mm];
\end{tikzpicture},
\,
then the four $\Aut(G)$-orbits are obvious. As minimal characteristic
subgroups are elementary abelian, the \emph{abelian} 4-orbit groups
are $(\Cyc_{p^3})^m$; $(\Cyc_{p^2})^k\times(\Cyc_p)^{m-k}$; and
$\Cyc_p^k\times\Cyc_r^{m-k}$ where $1\le k<m$, $p\ne r$ are prime and $m\ge1$.
The $\Aut(G)$-orbit lengths are obvious,
and the orders are $1,p,p^2,p^3$; $1,p,p,p^2$ and $1,p,r,pr$
respectively. 

By Lemma~\ref{L:C}, a nonabelian solvable 4-orbit group $G$ either has
four characteristic subgroups $G>M_1>M_2>1$  where $G'\in\{M_1,M_2\}$,
or is a UCS group (see~\cite{GPS}) with $G>G'=\Z(G)>1$.
\begin{table}[!ht]
\caption{Examples of non-nilpotent solvable 4-orbit groups where $p$ is a prime}
\label{T:4orbit}
\begin{center}
\begin{tabular}{clll}
  $G$&$\Aut(G)$-orbit lengths&Orders&Conditions\\
  \hline
  $(\Cyc_p)^d\rtimes\Cyc_r$&$1,p^d-1,ep^d,ep^d$&$1,p,r,r$&$e=\textup{ord}_p(r)=\frac{r-1}{2}\mid d$\\
  $(\Cyc_p)^d\rtimes\Cyc_r$&$1,2x,x^2,(r-1)p^d$&$1,p,p,r$&$\order_p(r)=\frac{r-1}{2}, x=p^{\frac{d}{2}}-1$\\
   $(\Cyc_{p^2})^d\rtimes\Cyc_r$&$1,p^d-1,p^{2d}-p^d,ep^{2d}$&$1,p,p^2,r$&$e=\order_p(r)=r-1\mid d$\\
   $(\Cyc_p)^d\rtimes\Cyc_{r^2}$&$1,p^d-1,(r-1)p^d,ep^d$&$1,p,r,r^2$&$e=\textup{ord}_p(r^2)=r(r-1)\mid d$\\
      &&&and $\textup{ord}_p(r)=r-1$\\
  $(\Cyc_p)^2\rtimes Q_8$&$1,p^2-1,p^2,6p^2$&$1,p,2,4$&$p\in\{3,5,7,11,23\}$
\end{tabular}
\end{center}
\end{table}
Hence a nilpotent 4-orbit group is a $p$-group of exponent
dividing $p^3$ with class at most~3.
We list in the next paragraph some infinite families of 4-orbit $p$-groups.
Some 4-orbit solvable non-$p$-groups are listed in Table~\ref{T:4orbit};
the first two lines are UCS groups. 
Dornhoff~\cite{Dornhoff73} studies
groups $N$ for which $\Aut(N)$ has a solvable subgroup, say $A$,
with four orbits on $N$
(i.e. three orbits on $N\setminus\{1\}$).
He lists $N$ in~\cite[Theorems~1.1,\,2.1]{Dornhoff73} and constrains the
structure of $N$ in~\cite[Theorems~3.1,\,4.1]{Dornhoff73}.
(No constraints are given in the case that $N$ has exponent $p$.)
The permutation groups $N\rtimes A\le\Sym(N)$ are
not always guaranteed to have rank 4. Hence obtaining a complete and
irredundant list of 4-orbit groups~$G$
may be quite difficult, especially in the case that each \emph{solvable}
subgroup of $\Aut(G)$ has more than 4 orbits on $G$.

There are many infinite families of 4-orbit $p$-groups. First, if $H$ is a
nonabelian 3-orbit $p$-group, and $E$ is an elementary abelian $p$-group,
then $G=H\times E$ is a 4-orbit group. Next, if $p$ is an odd prime
and $q=p^b$, then the group $G=\F_q^3\times\F_q^3\times\F_q$ in Example~\ref{L:GammaL3}
is a 4-orbit group of exponent $p$ with $\Aut(G)$-orbit lengths
$1,q-1,q^4-q,q^7-q^4$. Also, the extraspecial 2-groups $G_\eps$ of order
$2^{1+2k}$ with
$\eps=\pm$ and $(k,\eps)\ne(1,+)$ are 4-orbit groups of exponent $4$ by
Lemma~\ref{L1}. Then $|G_\eps|=2q^2$ where $q=2^k$ and
  the $\Aut(G_\eps)$-orbit
lengths are $1,1,q^2+\eps q-2,q(q-\eps)$ and the element orders $1,2,2,4$,
and $\Aut(G_\eps)$ is not solvable. Note that $q^2+\eps q-2=(q-\eps)(q+2\eps)$.

Some families of nonabelian solvable 4-orbit groups that are not $p$-groups
are listed in Table~\ref{T:4orbit}. The list uses, and extends~\cite[Theorem~4]{LM}.

  \section*{Acknowledgement}
  I am grateful to the referee for suggesting a number of improvements to this paper.


\end{document}